\RequirePackage{fix-cm}
\documentclass{svjour3}                     
\smartqed  
\usepackage{amssymb}
\usepackage{mathrsfs}
\usepackage{graphicx}
\usepackage{amsmath}
\usepackage{amsfonts}
\usepackage{bm}
\usepackage{subcaption}
\usepackage{hyperref}
\hypersetup{
    colorlinks=true,
    linkcolor=blue,
    filecolor=magenta,      
    urlcolor=cyan,
}
\usepackage{longtable,tabularx}
\usepackage{color}
\journalname{
Preprint}

\begin{document}

\title{Deep convolutional neural network for shape optimization using level-set approach}

\titlerunning{DL-ROM for shape optimization using level-set approach}        

\author{Wrik Mallik \and
        Neil Farvolden  \and
        Jasmin Jelovica \and
        Rajeev K. Jaiman
}


\institute{Wrik Mallik \\
           Jasmin Jelovica \\
           Rajeev K. Jaiman \\
            \at
               Department of Mechanical Engineering, \\
               The University of British Columbia, Vancouver, BC, Canada
           \and
           Neil Farvolden \\
           \at
              Institute for Aerospace Studies, \\
              University of Toronto, Toronto, ON, Canada 
}


\maketitle

\begin{abstract}
This article presents a reduced-order modeling methodology via deep convolutional neural networks (CNNs) for shape optimization applications. The CNN provides a nonlinear mapping between the shapes and their associated attributes while conserving the equivariance of these attributes to the shape translations. To implicitly represent complex shapes via a CNN-applicable Cartesian structured grid, a level-set method is employed. The CNN-based reduced-order model (ROM) is constructed in a completely data-driven manner thus well suited for non-intrusive applications. We demonstrate our ROM-based shape optimization framework on a gradient-based three-dimensional shape optimization problem to minimize the induced drag of a wing in low-fidelity potential flow. We show a good agreement between ROM-based optimal aerodynamic coefficients and their counterparts obtained via a potential flow solver. The predicted behavior of the optimized shape is consistent with theoretical predictions. We also present the learning mechanism of the deep CNN model in a physically interpretable manner. The CNN-ROM-based shape optimization algorithm exhibits significant computational efficiency compared to the full-order model-based online optimization applications. The proposed algorithm promises to develop a tractable DL-ROM driven framework for shape optimization and adaptive morphing structures.

\end{abstract}

\section{Introduction} 
Shape optimization has gained considerable interest over the last three-four decades with the advent of high-performance computing resources. The largely increased computing power has enabled the application of various numerical solution techniques like the finite element method (FEM), the finite difference method (FDM), etc., to explore solutions associated with very complex shapes. This has increased the recent popularity of shape optimization in structural optimization studies, whereas past studies relied solely on sizing optimization \cite{haftka1986structural}. Similarly, on the side of fluid mechanics, the high sensitivity of drag and aerodynamic performance on the airfoil shape has been known to the aerospace community for several decades \cite{jameson_1998,Giles_2000,pironneau_2004}. The advent of computational fluid dynamics (CFD) simulations has facilitated accurate computation of such sensitivities even for complex shapes. This has led to a growing interest in aerodynamic shape optimization over the last two decades \cite{reuther1996aerodynamic,samareh2004aerodynamic}.

Integrating high-fidelity computational structural mechanics (CSM) or CFD simulations within an optimization algorithm is often restricted by the huge computational expense of these simulations. This has restricted the computational model complexity that can be realistically employed in complex optimization studies \cite{mallik2015effect,de2021lightweight}. Some recent fluid dynamic shape optimization studies with high-fidelity CFD-based analyses can be found in Refs.~ \cite{garg2015high,lyu2015aerodynamic,he2019robust}. These studies have completely relied on a gradient-based search algorithm, which can very frequently converge to local optima for a non-convex, multi-objective shape optimization problem, especially when there are many design parameters \cite{chernukhin2013multimodality,streuber2020evaluating}. However, employing expensive metaheuristic algorithms that can convergence to the global optima remain largely non-tractable with high-fidelity solvers. This has motivated the present-day researchers in the development of computationally efficient reduced-order models (ROMs)/surrogate models, which can accurately emulate the computationally expensive full-order models (FOMs).

A variety of data-driven reduced-order modeling techniques exist in literature like Kriging \cite{krige1951statistical} (also known as Gaussian process regression) and
its variants \cite{schobi2015polynomial,bouhlel2016improving}, polynomial chaos expansions \cite{blatman2009adaptive}, and support vector machines \cite{feng2018multidisciplinary}. Some reviews have elaborated the applications of these methods for shape optimization applications \cite{koziel2013surrogate,han2015surrogate,yondo2018review}. Of all these methods, Kriging is the most widely used method in various engineering analysis and design tasks \cite{timme2011transonic,fan2019reliability,mallik2020kriging}. However, Kriging suffers from a variety of issues such as being poor at approximating discontinuous functions \cite{raissi2016deep}, difficulty in handling
high-dimensional problems \cite{perdikaris2017nonlinear}, costly to use in the presence of a large number of data samples \cite{keane2008engineering}, and being difficult
to implement for solving certain inverse problems with strong nonlinearities \cite{bonfiglio2018improving}.

With the advent of deep learning (DL), deep neural networks (DNNs) have been demonstrated as an effective approximation technique for most general complex functions \cite{cybenko1989approximation,hornik1991approximation}. They do not suffer from the many limitations of Kriging mentioned above and are replacing Kriging as the best ROM candidates for optimization applications \cite{bouhlel2020scalable,nagawkar2022aerodynamic}. Besides, increasingly wide feed-forward DNNs (also called multi-layer perceptrons (MLPs)), have been demonstrated to converge to Gaussian processes under broad conditions \cite{matthews2018gaussian}. This has led to the increased applications of MLPs \cite{kapania2018multidisciplinary,singh2021alga,wang2018database} and gradient-enhanced MLPs \cite{bouhlel2020scalable} for ROM-based shape optimization studies. The goal here is to construct an efficient DL-ROM based shape optimization of 3D wing with constraints. 

In this article, we present an efficient and unique approach for developing data-driven ROMs for shape optimization applications. We employ convolutional neural networks (CNNs) \cite{goodfellow2016deep} as a ROM for the higher-order computational model. CNNs have demonstrated excellent image learning capability in the commercial DL community because of the shared memory employed during the learning. The convolution (or cross-correlation) operation along the shared memory parameters endows CNNs with translational equivariance, whereby they can learn how the features of the input transform due to a translation of the input itself. Two or three-dimensional shapes can easily be picturized as pixels or voxels, respectively. Similarly, aerodynamic force coefficients, structural stresses, etc., associated with such shapes can be considered as output features of these pixels/voxels, and can thus be efficiently learned via CNNs. This is our motivation for employing CNNs as a ROM to learn the accurate mapping between the shapes and their associated features arising from various physical phenomena. A detailed discussion on how the convolution operation and translation equivariance forms strong geometric priors to enhance the learning capacity over MLPs is presented in Ref. \cite{goodfellow2016deep}. 

Level-set methods can provide an implicit representation of even complex shapes on a fixed grid. Often such an implicit representation leads to a coarse-grained and robust representation of shapes which would otherwise require a very fine mesh for explicit representation \cite{sethian1996theory}. This has often motivated the application of various forms of level-set representation approach for shape \cite{allaire2004structural,wang2007extended} and topology optimization \cite{wang2006radial,wang2007extended}. Although this remains one of the reasons for the application of level-set methods for the present article, the primary reason is somewhat different. Here we employ CNNs as our DL model. CNNs require input data in the form of a uniform point cloud on a Euclidean domain. Thus, the level-set representation of shapes on a relatively coarse, uniform Cartesian grid provides input information in a form directly employable by CNNs. Such features of the level-set representation for generating CNN input data have also been successfully demonstrated in \cite{miyanawala2018novel}. There the CNN was employed for learning the aerodynamic force coefficients obtained from Navier-Stokes simulations of two-dimensional flow around bluff-bodies \cite{miyanawala2018novel}. Here we employ it for the first time for implicit three-dimensional shape representation on a CNN-compatible grid to develop a ROM for shape optimization applications.

In the present study, we specifically employ deep CNNs with level-set shape representations to obtain a mapping between three-dimensional wing configurations and their aerodynamic force coefficients. The FOM employed here is a low-fidelity potential flow solver, which can predict the lift and induced drag generated by the wing at low speeds. However, it could easily be replaced by higher-fidelity CFD simulations or experiments without any loss of generality, as the CNN-based mapping is completely data-driven. The FOM is employed for generating an adequate number of input and output data sets to train the deep CNN model. The CNN-based ROM is subsequently employed for non-intrusive shape optimization of a wing shape to minimize the wing's drag coefficient, while prescribed constraints are applied on the wing lift coefficient and minimum thickness ratio.

The methodology of developing the CNN-based ROM with level-set methods, its application in the constrained optimization study, and the accuracy of the results predicted by the CNN model are discussed in the paper. We specifically demonstrate the optimization applications with gradient-based methods but the non-intrusive application enables seamless integration with other gradient-free optimization algorithms like trust region-based methods \cite{marques2021non} and metaheuristic algorithms \cite{mallik2015effect,de2019structural,samanipour2020adaptive,de2021lightweight}. We also explore the feature space of various convolution layers of the deep CNN model to provide interpretability to the deep CNN's learning mechanism. The computational efficiency of the present approach compared to FOM-based shape optimization is also discussed. Such application of deep CNNs with level-set methods as ROMs for shape optimization is presented here for the first time.

\section{Uniqueness of the proposed reduced-order modeling approach}
There has been a recent application of NNs as ROMs for shape optimization applications \cite{bouhlel2020scalable,nagawkar2022aerodynamic}. These articles have focused on MLPs as their DL framework and have specifically employed a gradient-enhanced version of such networks. While such methods have shown good accuracy, they have been explicitly developed for gradient-based optimization algorithms, which may not be adequate in many cases \cite{chernukhin2013multimodality,streuber2020evaluating}. Another set of data-driven ROMs employed recently was based on proper orthogonal decomposition (POD) bases \cite{yao2020reduced} and hyper-reduction based on POD modes \cite{marques2021non}. While the latter method was developed for a non-intrusive application, linear projection methods like POD are less efficient for certain convection-dominated fluid mechanics problems \cite{mojgani2021low}. 

Our proposed model order reduction approach is unique to such recent articles in the following ways:
\begin{itemize}
    \item We employ CNNs as our DL model, which is significantly more efficient than linear projection-based model reduction approaches for even the most complex physical cases \cite{mojgani2021low}.
    \item We employ CNNs to obtain a map between shapes represented on a uniform Cartesian grid and their associated attributes. CNNs are specifically suited for such learning tasks aided by certain geometric priors. These are discussed later in the article.
    \item We employ level-set methods for implicit shape representation, which allows us to represent shape information on a CNN-employable grid. Furthermore, such implicit representation can be obtained on a coarser mesh with fewer parameters compared to explicit parametric representation with free-form-deformation parameters employed for three-dimensional shape optimization study \cite{samareh2004aerodynamic,garg2015high}.
\end{itemize}

\section{Deep learning-based model order reduction methodology}
The DL-based model order reduction methodology consists of two parts, the implicit shape representation, and the development of the DL model for learning the functional relation between shapes and their associated attributes. These are explained below in detail.

\subsection{CNN: a geometric DL model}
In the present study, we want our DL model to successfully learn the changes in two or three-dimensional shape attributes with changes in the shapes. This requires the DL model to satisfy certain prior geometric assumptions. CNNs can satisfy these geometric priors, resulting in their tremendous success for similar applications in the field of computer vision \cite{lecun2015deep,saha2022practical}. These geometric priors required for the supervised learning task are discussed below.

Let us consider a compact $d$-dimensional Euclidean domain $\Omega = {[0,1]}^d \subset \mathbb{R}^d$ on which square integrable functions $f \in L^2 \left(\Omega \right)$ are defined. We consider a generic supervised learning environment for obtaining an unknown function $y: L^2 \left(\Omega\right) \rightarrow \mathcal{Y}$ on a training set,
\begin{equation}
    {\{f_i \in L^2 \left(\Omega\right),y_i=y\left(f_i\right) \}}_{i \in \mathcal{I}}.
\end{equation}
We are interested in regression tasks, thus $\mathcal{Y}=\mathbb{R}^m$. For the proposed learning problem, the unknown function $y$ needs to satisfy the crucial prior assumptions explained below.

\emph{Stationarity}: Our first assumption is that depending on the task, the function $y$ is either invariant or equivariant with respect to translations. Here
\begin{equation}
    \mathcal{T}_v f(x) = f(\left(x-v\right), \quad x, v, \in \Omega,
\end{equation}
is defined as a translation operator acting on functions $f \in L^2 \left(\Omega\right)$. When we consider translation invariance, we have $y(\mathcal{T}_v f) = y(f)$ for any $f \in L^2 \left(\Omega\right)$ and $v \in \Omega$. For translational equivariance, we have $y(\mathcal{T}_v f) = \mathcal{T}_v y(f)$.

\emph{Local deformations and scale separation}: We also desire $y$ to be stable to local deformations $\mathcal{L}_\tau$ , where $\tau: \Omega \rightarrow \Omega$ is a smooth vector field, acts on $L^2 \left(\Omega\right)$ as $\mathcal{L}_\tau f(x) - f\left(x - \tau(x)\right)$.

We briefly explain below how the CNN satisfies such prior assumptions during the forward propagation phase of the supervised learning algorithm. The supervised learning algorithm also includes a backpropagation phases, which ensures the CNN framework learns optimal parameters to minimize a prescribed loss function.

\subsubsection{Forward propagation phase}
In the present application, stationarity and stability to local translations are both leveraged in the CNN framework as follows. The CNN network consists of several convolutional layers of the form $C(\bm{f})$ which sequentially performs a set of convolution operation of the form $\bm{g}=\kappa_\Lambda(\bm{f})$, a batch-normalization operation $\beta$, and a point-wise non-linearity $\xi$, acting on a $p$-dimensional input $\bm{f}(x) = \left( f_1(x), \ldots , f_p(x) \right)$,
\begin{equation}
    C(\bm{f}) = \xi\left(\beta\left(\kappa_\Lambda(\bm{f})\right)\right).
\end{equation}
The convolution operation $\kappa_\Gamma(\bm{f})$ operates by applying a set of kernels (or filters) $\Lambda = \left(\lambda_{l,l^{'}}\right), l=1, \ldots, ,p, l^{'} = 1, \ldots, p$, 
\begin{equation}
    g_l(x) =  \sum_{l^{'}=1}^p \left(f_{l^{'}}*\lambda_{l,l^{'}}\right)(x),
\end{equation}
producing a $q$-dimensional output $\bm{g}(x) = \left(g_1(x), \ldots, g_q(x)\right)$, often referred to as the feature maps. Here,
\begin{equation}
    \left(f_{l^{'}}*\lambda\right)(x) = \int_{\Omega} f\left(x-x{'}\right)\lambda(\left(x{'}\right)dx{'},
\end{equation}
denotes a standard convolution operation.

Batch normalization is an operation used to stabilize the training process and ensure faster training convergence. This is achieved by normalizing the output $\bm{g}$ of a convolution operation by their mean $\mu$ and variance $\sigma$,
\begin{equation}
    \beta(\bm{g}) = \frac{\bm{g} - \mu}{\sqrt{\sigma} + \epsilon},
\end{equation}
where $\epsilon$ is a small positive value used to improve numerical stability when $\sigma \rightarrow 0$.
Further details of the batch normalization process can be obtained in Ref. \cite{ioffe2015batch}. Finally, we introduce nonlinearity in our learning mechanism with a leaky rectified linear unit (leakyReLU) \cite{maas2013rectifier}.  This is performed via the nonlinear activation layer as follows,
\begin{equation}
    \xi(x) = \max\left(x, \gamma x\right),
\end{equation}
where $\gamma$ is a scaling factor. 

A CNN model composed of $K$ convolutional layers put together has a general hierarchical representation,
\begin{equation}
    \label{eq:CNN_model}
    U_\Theta(\bm{f}) = \left(C_{(K)} \circ \ldots C_{(2)} \circ C_{(1)}\right),
\end{equation}
where $\Theta = \{\Lambda^{(1)}, \ldots,  \Lambda^{(K)}\}$ is the set of all network parameters (all the kernel coefficients). The model is considered deep when it consists of multiple convolutional layers. For the present purpose we assume a CNN model with four convolutional layers to be deep. Further details about the CNNs and their geometric learning capability on an Euclidean domain can be obtained in Ref.~\cite{bronstein2017geometric}.

\begin{figure}[ht!]
\centering
\includegraphics[width=0.98\linewidth]{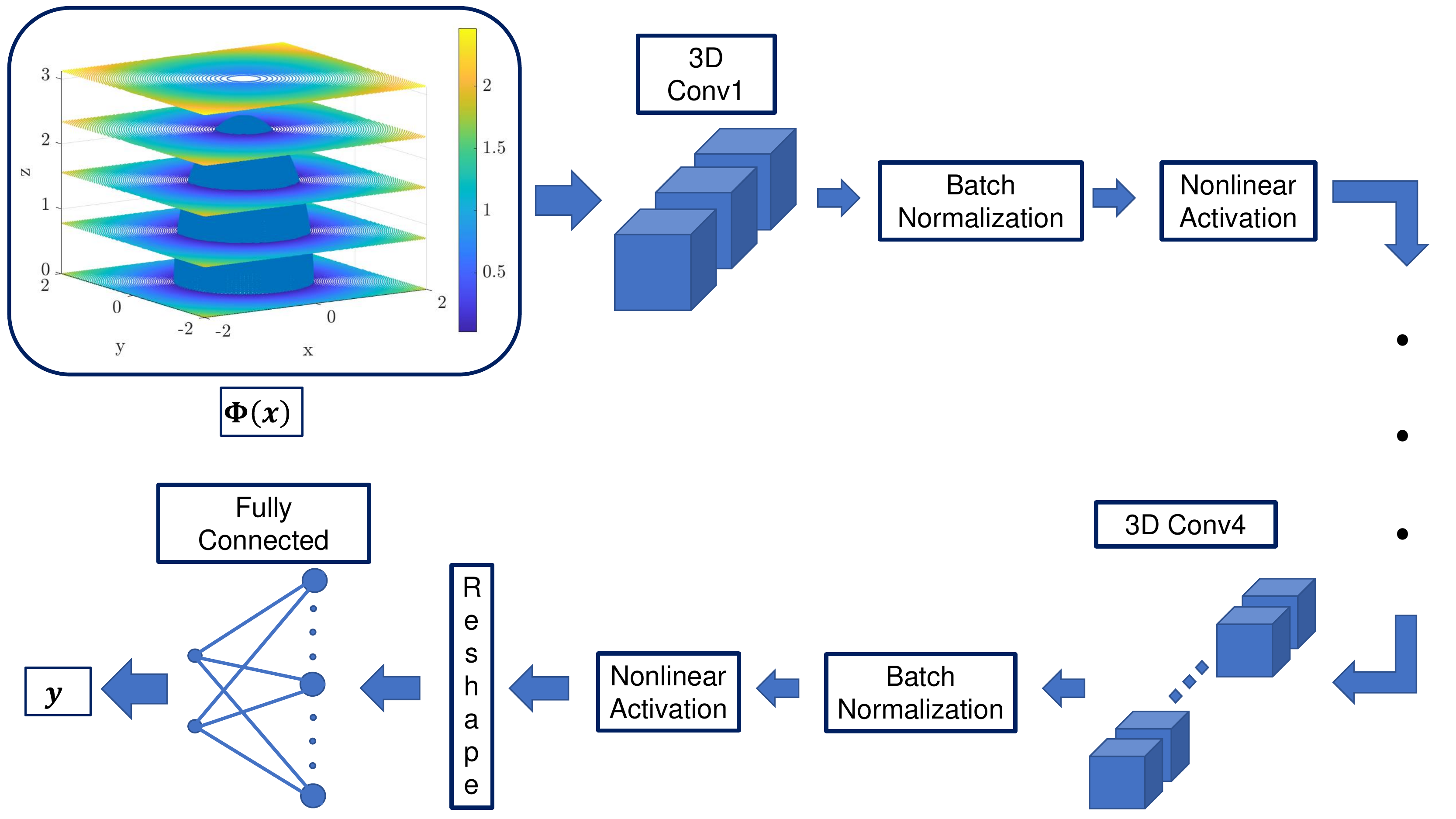}
\caption{Deep CNN architecture for shape optimization}
\label{fig:CNN_architecture}
\end{figure}

The deep CNN framework employed in this study consisted of four convolutional layers is illustrated in Fig.~\ref{fig:CNN_architecture}. The input data is the level-set representation, $\mathrm{\Phi}$, of the distance from the shape boundary. It is discussed next. The convolutional layers are followed by a fully connected layer to ensure a consistent dimension of the CNN output.

\subsubsection{Backpropagation phase}
The role of backpropagation is to compute the errors between the target, $\bm{y}$, and the predicted output, $\hat{\bm{y}}$, and update the trainable network parameters to minimize the errors. The target $\bm{y}$ is obtained from the FOM for the given inputs. The errors are represented via a loss function computing mean square error over a batch size $n_b$,
\begin{equation}
    \mathscr{L}\left(\hat{\bm{y}},\bm{y}\right) = \frac{1}{2} \sum_{i=1}^{n_b} {\left(\hat{\bm{y}}_i - \bm{y}_i \right)}^2.
\end{equation}
The weights and biases of the convolutional kernels and the fully connected layer are denoted generally as $W$ and $b$ from the trainable network parameters. These are updated iteratively via the gradient descent-based optimization Adam \cite{kingma2014adam}. Adam is run over several iterations to achieve convergence. To accelerate the convergence, the input dataset is divided into several mini-batches, and the loss function is evaluated for each mini-batch. The order of the mini-batches is selected randomly leading to a stochastic optimization process.

The weight (and similarly the bias) update process is shown below
\begin{equation}
    W^{t+1} = W^{t} - \frac{\alpha u^t}{\sqrt{v^t} + \epsilon}, 
\end{equation}
where $\alpha$ is the initial learning rate, $u_w^t$ is an adaptive gradient of the loss with respect to the weight, and $v_w^t$ is a moving average of the element-wise squares of the loss gradients. The superscript $t$ denotes the iteration of the optimization process. Here, we consider $\alpha$ as a hyperparameter that can be tuned to improve the network's predictions. $u_w^t$ and $v_w^t$ are obtained as,
\begin{equation}
    \begin{split}
        u_w^t &= \beta_1 u^{t-1} + (1-\beta_1)\frac{\partial \mathscr{L}}{\partial W} \\
        v_w^t &= \beta_2 v^{t-1} + (1 - \beta_2) {\left(\frac{\partial \mathscr{L}}{\partial W}\right)}^2.
    \end{split}
\end{equation}
The adaptive gradient $u_w^t$ helps to reduce the oscillations of the stochastic optimizer as it converges to an optimum by reducing the step size along the path of the steepest gradient. $v_w^t$ serves a similar purpose but instead adaptively decays the learning rate of the weight update with increasing iterations, as the optimizer converges to the optimal weights. $\beta_1$ and $\beta_2$ are the decay factor of the gradients and squared gradients of the loss and epsilon is a stabilizer of the numerical process for very small values of $v_w^t$.

\subsection{Level-set method for shape representation}
Here we employ a level-set method to represent the three-dimensional shape on a structured three-dimensional Cartesian grid. Level-set methods are well-known for their capacity to implicitly represent even complex shapes on a coarser mesh \cite{sethian1996theory} compared to explicit techniques like free-form deformation variables. The computational efficiency of such coarse-grained, implicit representation has motivated their extensive application for shape and topology optimization studies \cite{allaire2004structural,wang2006radial,wang2007extended}. 

\begin{figure}[ht!]
\centering
\subfloat[Ellipsoid shape]{\includegraphics[width=0.47\linewidth]{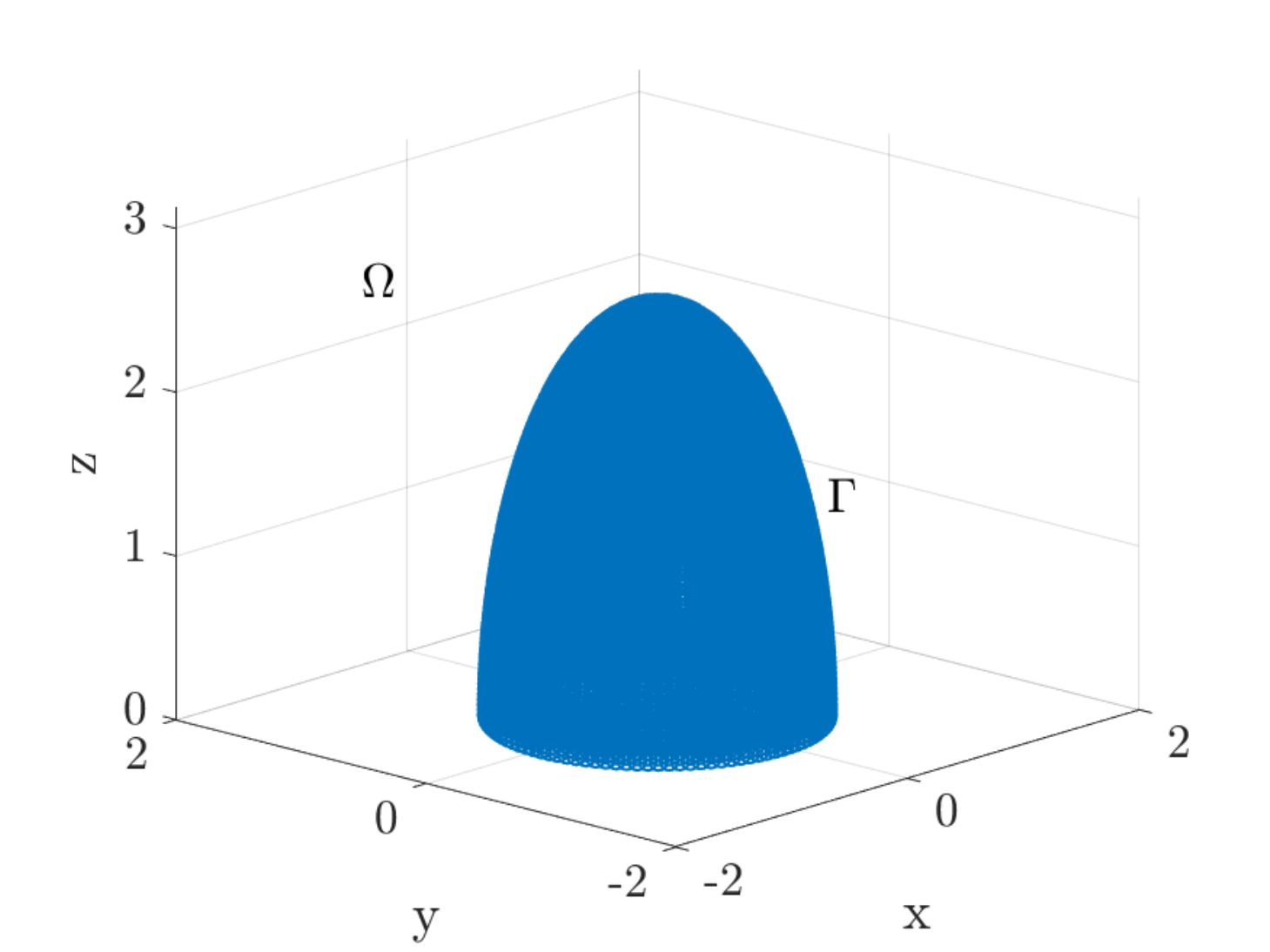}}
\subfloat[Level-set representation]{\includegraphics[width=0.51\linewidth]{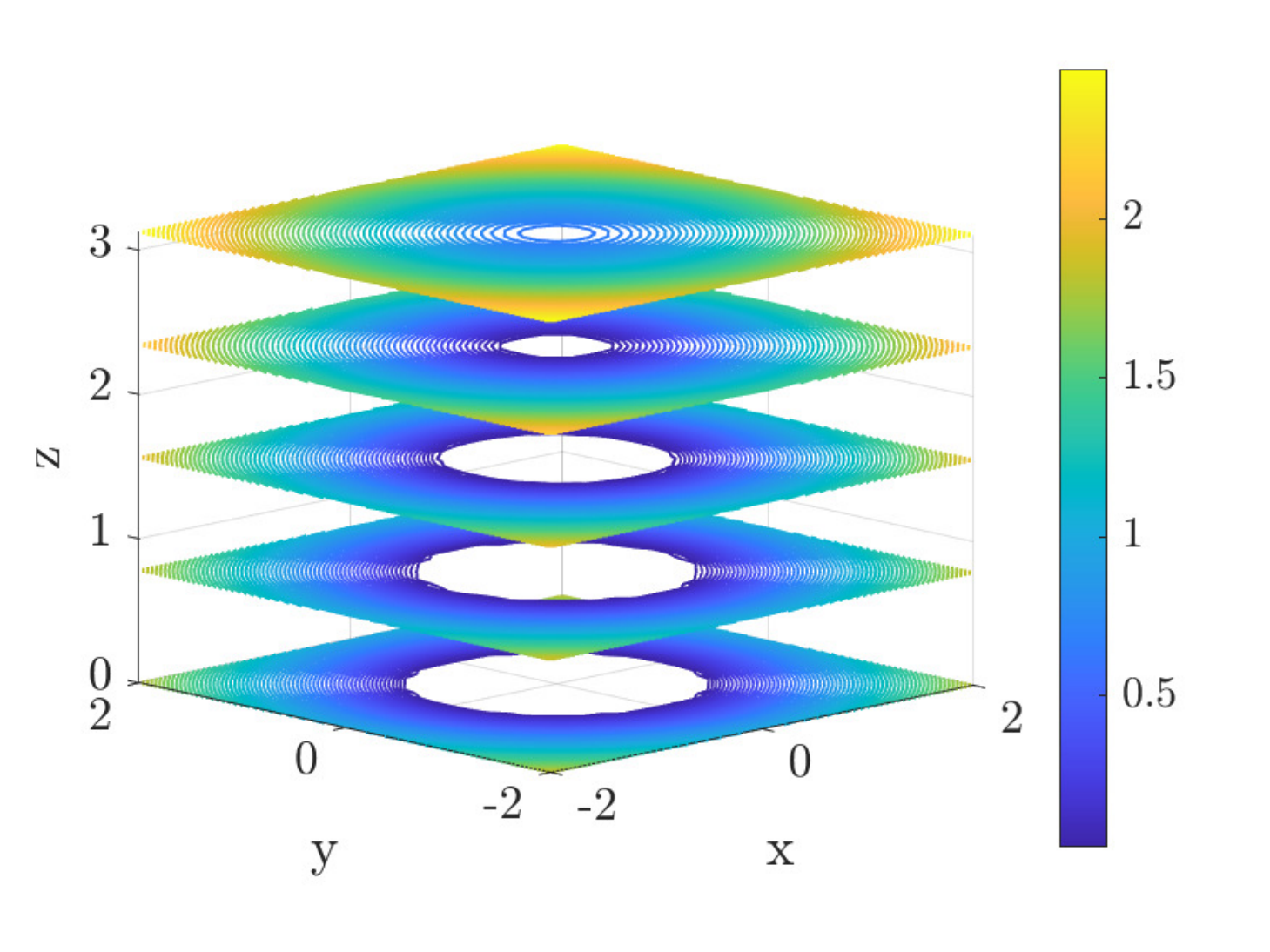}}
\caption{Level-set implicit representation of an arbitrary ellipsoid shape}
\label{fig:ellipsoid_level_set}
\end{figure}

We represent the three-dimensional shape via level-sets of the minimum Euclidean distance of any point on the level-set grid from the shape boundary. The boundary and domain inside it are represented by zero level-sets. These are formulated as follows:
\begin{equation}
    \label{eq:level-set_def}
    \mathrm{\Phi}\left(\bm{x};\Gamma\right) = \begin{cases}
    \min \left({\Vert \bm{x} - \bm{x}_\Gamma \Vert}_2\right), &\quad \forall \bm{x} \in \Omega\\
    0 &\quad \forall \bm{x} \in \Gamma \\
    0 &\quad \forall \bm{x} \in D \\
    \end{cases}
    ,
\end{equation}
where $\bm{x}$ represents the Cartesian coordinates, $\Gamma$ represents the shape boundary, $\Omega$ represents the domain outside the boundary, and $D$ represents the domain inside the boundary $\Gamma$. We show the shape boundary and the domain outside the boundary for an arbitrary ellipsoid configuration in Fig.~\ref{fig:ellipsoid_level_set} (a). The slices of the various level-set iso-surfaces for the ellipsoid configuration are presented at various locations along the $z$-axis in Fig.~\ref{fig:ellipsoid_level_set} (b). 

Apart from the implicit shape representation, the level-set method serves another purpose for the present application. Here we intend to employ CNNs as our DL model for the supervised learning of how physical attributes associated with various shapes change with their geometry. The CNNs learn such maps on a Euclidean domain and require a point cloud of the input data on a uniform, structured grid. The level sets of distance function obtained on a uniform Cartesian grid thus provide our shape input in a form directly employable by CNNs.

\section{DL-ROM-based shape optimization methodology}
\begin{figure}[ht!]
\centering
\includegraphics[width=0.98\linewidth]{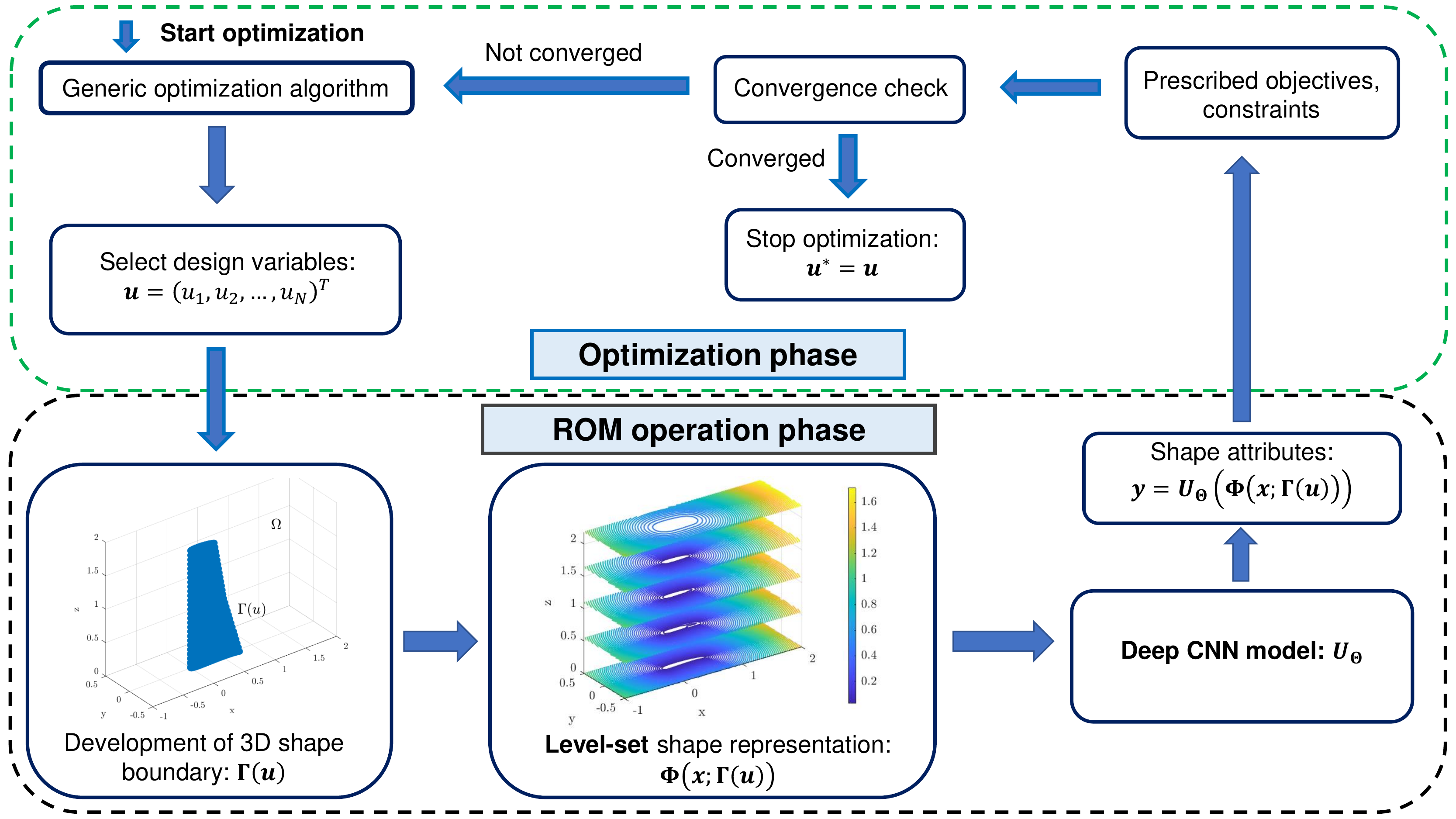}
\caption{DL-ROM-based shape optimization framework}
\label{fig:DL-ROM-based_opt}
\end{figure}

The DL-ROM-based shape optimization methodology can be divided into two phases. One is an optimization phase, which involves defining and updating the design variables, and steering the search via gradients of the objective and constraints. The second phase is the DL-ROM operation phase where the level-set shape representation and CNN predictions are carried out as a ``black-box" with respect to the optimization algorithm. Fig.~\ref{fig:DL-ROM-based_opt} shows a schematic demonstrating the flow of information in the shape optimization process. The deep CNN model block in the ROM operation phase of Fig.~\ref{fig:DL-ROM-based_opt} is a compressed representation of the CNN framework presented in Fig.~\ref{fig:CNN_architecture}. Interchange of information takes place at two stages, once when the design variables selected by the optimization algorithm are passed to the DL-ROM, and again when the shape attributes associated with these design variables are passed back to the optimization phase by the DL-ROM model. The DL-ROM is thus employed in a non-intrusive manner. A detailed explanation of the two phases of the complete DL-ROM-based shape optimization methodology is provided below.

\subsection{Non-intrusive optimization phase}
The shape optimization problem is formulated as,
\begin{equation}
    \begin{split}
        \label{eq:gen_opt_formulation} 
        \underset{u \in \mathscr{V}}{\text{minimize}} & \quad \mathcal{J} \left(\bm{y}\left(\bm{u}\right), \bm{u}\right)\\
        \text{subject to} & \quad
        \mathcal{C_E}\left(\bm{y}\left(\bm{u}\right), \bm{u}\right) = 0,\\
        & \quad \mathcal{C_I}\left(\bm{y}\left(\bm{u}\right), \bm{u}\right) \leq 0.
    \end{split}
\end{equation}
Here, $\bm{u}$ are a set of $N$ design variables, $\bm{u}=\left(u_1, u_2, \dots, u_N\right), \bm{u} \in \mathbb{R}$, which can change the three-dimensional shapes. $\bm{y}\left(\bm{u}\right)$ are the attributes of the shapes thus created and are obtained directly from $\bm{u}$ via the black-box DL-ROM (Eq.~\ref{eq:black-box_surrogate}) during the DL-ROM operation phase, $\mathcal{J}$ is the objective function, $\mathcal{C_E}$ is the equality constraint and $\mathcal{C_I}$ is the inequality constraint. $\mathscr{V} \subset \mathbb{R}^{N}$ represents the $N$-dimensional parameter space of all possible combination of design variables.

Eligible solutions to this problem are the ones that satisfy all constraints and are members of feasible set $\bm{U}$
\begin{equation}
    \label{eq:opt_bounds}
    \bm{U}=\{\wp \left({\left[u_1,\ldots,u_N\right]}^T\right)|\bm{b_l}\leq\bm{u}\leq\bm{b_u}\},
\end{equation}
which contains all variable permutations $\wp$ of $\bm{x}$ between their lower and upper bounds, $\bm{b_l}$ and $\bm{b_u}$, respectively. Thus, it contains all possible design alternatives. 

The DL-ROM-based shape optimization methodology presented above is a generic formulation for non-intrusive DL-ROM applications, compatible with any generic optimization algorithm. In the present study, we employ a gradient-based algorithm to obtain our optimal configuration. The gradient of the objective function is computed as
\begin{equation}
    \frac{\mathrm{d} \mathcal{J}}{\mathrm{d} \bm{u}} = \frac{\mathrm{\partial} \mathcal{J}}{\mathrm{\partial} \bm{u}} + \frac{\mathrm{\partial} \mathcal{J}}{\mathrm{\partial} \bm{y}} \frac{\mathrm{d} \bm{y}}{\mathrm{d} \bm{u}},
\end{equation}
where $\frac{\mathrm{\partial} \mathcal{J}}{\mathrm{\partial} \bm{u}}$ and $\frac{\mathrm{d} \bm{y}}{\mathrm{d} \bm{u}}$ are computed directly via any routine numerical differentiation scheme. The shape sensitivities of the constraints, $\frac{\mathrm{d} \mathcal{C_E}}{\mathrm{d} \bm{u}}$ and $\frac{\mathrm{d} \mathcal{C_I}}{\mathrm{d} \bm{u}}$, are computed numerically analogous to the shape sensitivity of the objective function. Such computation of sensitivities are termed as monolithic differentiation in the gradient-based shape optimization literature \cite{martins2013review}, and is employed when the computational model is treated as a "black-box" model.

\subsection{DL-ROM operation phase}
The DL-ROM operation phase in shown schematically in Fig.~\ref{fig:DL-ROM-based_opt}. Mathematically the DL-ROM operation during an online application can be formulated as follows. We have defined a set of design variables $\bm{u}$ during the optimization phase used for creating a new three-dimensional configuration, and consequently change the boundary $\Gamma$ of the shape. Thus, for optimization purpose, we can represent the boundary as:
\begin{equation}
    \label{eq:opt_shape_boundary}
    \Gamma \equiv \Gamma(\bm{u}).
\end{equation}
For optimization application, we can substitute Eq.~\ref{eq:opt_shape_boundary} in Eq.~\ref{eq:level-set_def}  to obtain a design variable-based level-set representation:
\begin{equation}
    \label{eq:level-set_def_opt}
    \mathrm{\Phi}\left(\bm{x};\Gamma (\bm{u}) \right) = \begin{cases}
    \min \left({\Vert \bm{x} - \bm{x}_{\Gamma(\bm{u})} \Vert}_2\right), &\quad \forall \bm{x} \in \Omega\\
    0, &\quad \forall \bm{x} \in \Gamma(\bm{u}) \\
    0, &\quad \forall \bm{x} \in D \\
    \end{cases}
    .
\end{equation}
We subsequently provide $\Phi\left(\bm{x};\Gamma (\bm{u}) \right)$ as inputs to the CNN model (Eq.~\ref{eq:CNN_model}) to obtain the shape attributes $\bm{y}$,
\begin{equation}
    \label{eq:opt_surrogate_output}
    \bm{y} = U_\Theta\left(\mathrm{\Phi}\left(\bm{x};\Gamma (\bm{u}) \right)\right).
\end{equation}
For non-intrusive application in an optimization algorithm, the DL-ROM presented in Eq.~\ref{eq:opt_surrogate_output} can be substituted by an abstract form:
\begin{equation}
    \label{eq:black-box_surrogate}
    \bm{y} = \mathcal{F}\left(\bm{u}\right).
\end{equation}
These shape attributes $\bm{y}$ are now fed back into the optimization phase for computation of the objective function and constraint values, and their sensitivities with respect to the design variables $\bm{u}$ if a gradient-based optimization algorithm is employed.

\section{Test problem}
Here we apply our CNN-based DL-ROM for the shape optimization of a three-dimensional wing. Our objective is to minimize the induced drag of the wing $C_{D_i}\left(u\right)$, subject to certain constraints, all of which are functions of the design variables $u$. To perform a non-intrusive gradient-based shape optimization, we first begin with our offline training phase where the CNN is trained with an adequate number of wing shapes and the associated lift and drag coefficients, as the input and output data sets, respectively. Here, the level-set methods are used to obtain the CNN input data for the wing shape on a structured Cartesian grid. A potential flow solver is employed to obtain the lift and drag coefficients for the wing configurations. After the training, we perform the DL-ROM-based non-intrusive optimization with a gradient-based algorithm according to the methodology as outlined in Fig.~\ref{fig:DL-ROM-based_opt}. Details of each of the various operations carried out at various phases of the whole methodology are provided as follows.

\subsection{Parameterization of aerodynamic shape}
In the current study, the wing shape is obtained by linearly interpolating the airfoils at the four sections along the span. These include the wing root and the tip sections and two other equally spaced intermediate sections. The airfoil shapes are formed by combining two Ferguson splines, one each for the top and bottom surfaces \cite{sobester2008quest,sobester2014aircraft}. Ferguson airfoils facilitate a relatively simple parameterization of the wing shape and yet enable the generation of a range of subsonic airfoils. Each of the wing sectional airfoils can be defined by eight shape variables. Amongst these, six variables satisfy the tangency and airfoil wedge angle constraints \cite{sobester2008quest,sobester2014aircraft}, while the seventh and eighth variables define the sectional twist and airfoil chord, respectively. Thus, the four sections combine to form thirty-two design variables. The thirty-third and final design variable of the design variable set $\bm{u}$ is the angle of incidence of the whole wing.

\subsection{Aerodynamic analysis}
The aerodynamic analysis for all three-dimensional wing configurations was conducted using linear lifting line theory (LLT) to provide a relatively quick generation of training data. We will consider low-speed flow, attached flow conditions, and unswept, large aspect-ratio wing configurations. Under such conditions, lifting line theory can provide a fairly accurate solution for the induced drag, $C_{D_i}$. Furthermore, most wing configurations selected here will be relatively thin. Under such conditions, the induced drag is the major contributor to the total drag on a wing. Thus, LLT can be considered an acceptable numerical solver to obtain the CNN training data. The lift and induced drag coefficients were computed using 20 points evaluated along the half span and 20 even-numbered Fourier coefficients. The sectional aerodynamic properties used in the LLT analysis were acquired via panel method solutions. 200 panels along the combined upper and lower surfaces provided reasonable convergence of the sectional potential flow solutions. The theory of the well-known LLT and panel method techniques are not provided here for brevity. Detailed information can be found in Ref.~\cite{anderson2010fundamentals}. 

It is important to note that the DL-ROM is completely data-driven. Thus, the CNN application methodology will not be affected if higher fidelity CFD or experimental data were augmented to the training data set. Thus, we can use the present data generation approach to demonstrate the model reduction capabilities of the CNN without any loss of generality.

\subsection{Level-set wing representation as CNN input}
As explained earlier the level-sets of the distance function are computed on a three-dimensional Cartesian grid to be provided as input data to the CNN. Thus, the level-sets provide a discrete three-dimensional representation $\bm{\Phi} \in \mathbb{R}^{m\times n \times p}$, where $m$, $n$ and $p$ represent the discretization dimensions in the $x$, $y$ and $z$ directions. As shown earlier in Eq.~\ref{eq:level-set_def_opt}, the level-sets provide the Euclidean distance to the level-set grid point from the nearest point on the wing boundary, for all points external to the boundary.

The discrete wing distance level-set are computed on the bounded domain $\Omega$, 
\begin{equation}
    \Omega(x,y,z) \in \{[-1,3]\times[0,2]\times[-0.5,0.5]\},
\end{equation}
along the $x$, $y$, and $z$ directions, respectively. Following aerodynamic conventions, these represent the chord-wise, span-wise, and transverse directions. For the present test problem, we have $m=48$, $n=16$ and $p=32$. Higher resolutions are chosen for the $x$ and $z$ directions because of the more pronounced variability in the geometry along these axes. In the $y$ direction, most of the axial variability in geometry comes from gradual changes in thickness and chord, and so the computational cost related to input size may be reduced without significant losses in model fidelity. In addition, the wing surface $\Gamma$ is discretized using 50 chord-wise reference points (for each the top and bottom surface) and 50 span-wise reference points for computation of the Euclidean distance from the surface. This wing surface definition is independent of the wing level-set grid resolution.

\subsection{Wing shape optimization problem}
Our objective is to minimize the induced drag coefficient $C_{D_i}$ subject to constraints on the design lift coefficient $C_L$, and the wing's minimum thickness to chord ratio, $t_c$, by changing the wing's shape design variables $\bm{u}$. Thus, referring to the general non-intrusive shape optimization formulation presented in Eq.~\ref{eq:gen_opt_formulation}, for the present case
\begin{equation}
    \begin{split}
        \label{eq:opt_test_problem} 
        \mathcal{J} &=  C_{D_i}(\bm{u}) \\
        \mathcal{C_I}(1) &= {t_c}^*- \max \left(t_c(\bm{u})\right)\\
        \mathcal{C_I}(2) &= \lvert{C_L}^* - C_L(\bm{u})\rvert - \epsilon.
    \end{split}
\end{equation}
It is worth mentioning that there are no equality constraints for the present case of wing shape optimization.

The first constraint is applied on thickness to chord ratio $t_c$ to ensure that its maximum value amongst all sections of the wing satisfies the minimum thickness to chord ratio, ${t_c}^*$. It is employed to obtain realistic manufacturing sizes and also indirectly ensure that the wing sections will be able to satisfy the stress requirements of a post-optimization structural analysis. The second constraint is applied on the lift coefficient $C_L$ to ensure that it lies within a tolerance level, $\epsilon$, of the design lift-coefficient $C_L^*$. This ensures that the designs meet the required level-flight design lift conditions. For our study, the minimum thickness to chord ratio ${t_c}^*$ is selected as 8\% and the design lift-coefficient ${C_L}^*=0.42$. We generally select the tolerance $\epsilon=0.05 C_L$.

\section{Results}
We first obtain our training data sets via FOM computations. These data are used for training, validation and testing our deep CNN model. We then analyze the feature space of the various composite convolutional layers of our deep CNN framework to provide interpretability to the learning mechanism of CNN. We subsequently employ our CNN-based DL-ROM with the level-set shape representation for online applications with a gradient-based optimization algorithm. Each of these processes is detailed below.

\subsection{CNN training data generation}
To obtain the CNN training data sets, the level-set of the inputs must be first generated for a wide range of wing shape variables. The bounds for Ferguson airfoil geometry variables were selected as per the recommendations of Sobester and Forrester \cite{sobester2014aircraft}. The lower and upper bounds of the wing's incident angle of attack were selected as -3 and 10 degrees, respectively. This ensures a sufficiently wide range of lift and coefficient values for the CNN to accurately learn the forward map from the geometry parameters to the geometry features.

A Latin Hypercube (LH) sampling of the thirty-three dimension design variable space was performed to generate 17,500 wing configurations. Amongst them 15,000 were selected for the training set, 1,500 designs were selected for the validation set, and the remaining 1,000 designs were selected for the test set. The training, validation and test sets were also selected randomly from the total data sets. Furthermore, the mean and variance of each individual data set were inspected to ensure that they represent similar random distributions. The LLT was employed to obtain the $C_L$ and $C_{D_i}$ for all these designs and the distance function level sets were also computed. The distance function level-sets and the aerodynamic force coefficients are the input and output data sets for the CNN, respectively.

\subsection{Training, validation and testing of the CNN}
\begin{figure}[hbt!]
\centering
\includegraphics[width=0.98\linewidth]{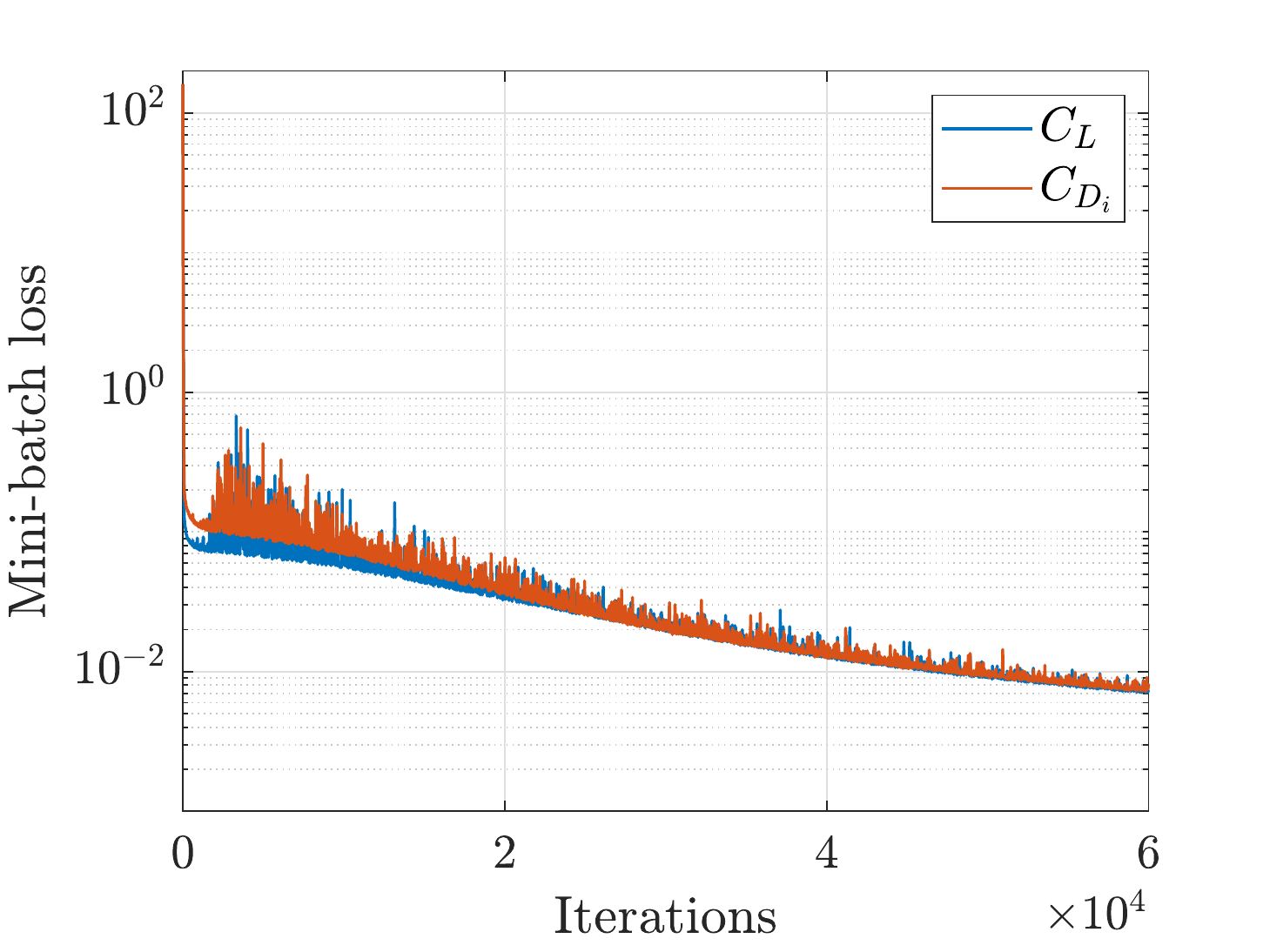}
\caption{Training convergence of CNN: $C_L$ and $C_{D_i}$}
\label{fig:train_conv}
\end{figure}

The CNN was trained on a single NVIDIA V100 Tensor Core GPU. Matlab's deep learning toolbox \cite{matlabdl} was employed to develop and train the CNN architecture. Separate CNN models were employed for predicting $C_L$ and $C_{D_i}$. The training convergence for both the $C_L$ and $C_{D_i}$ CNNs are presented in Fig.~\ref{fig:train_conv}. It is observed that over four orders of mini-batch loss could be reduced by training the CNN for 60,000 iterations. The figure also shows that after about 50,000 iterations the reduction of the mini-batch loss is gradual. Thus, satisfactory  convergence of the CNN model can be expected after 60,000 iterations.

\begin{figure}[hbt!]
\centering
\subfloat[$C_L$]{\includegraphics[width=0.49\linewidth]{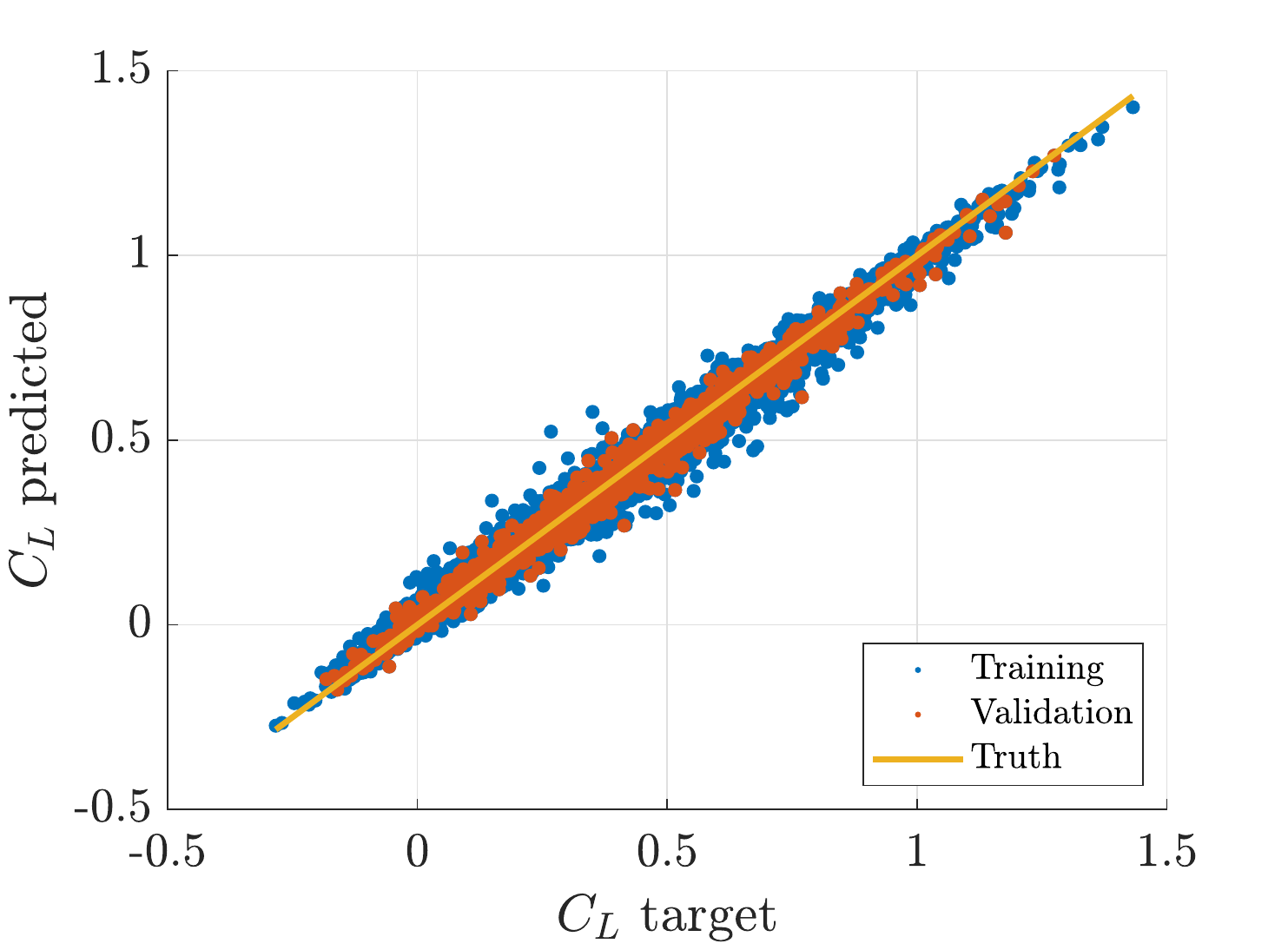}}
\subfloat[$C_{D_i}$]{\includegraphics[width=0.49\linewidth]{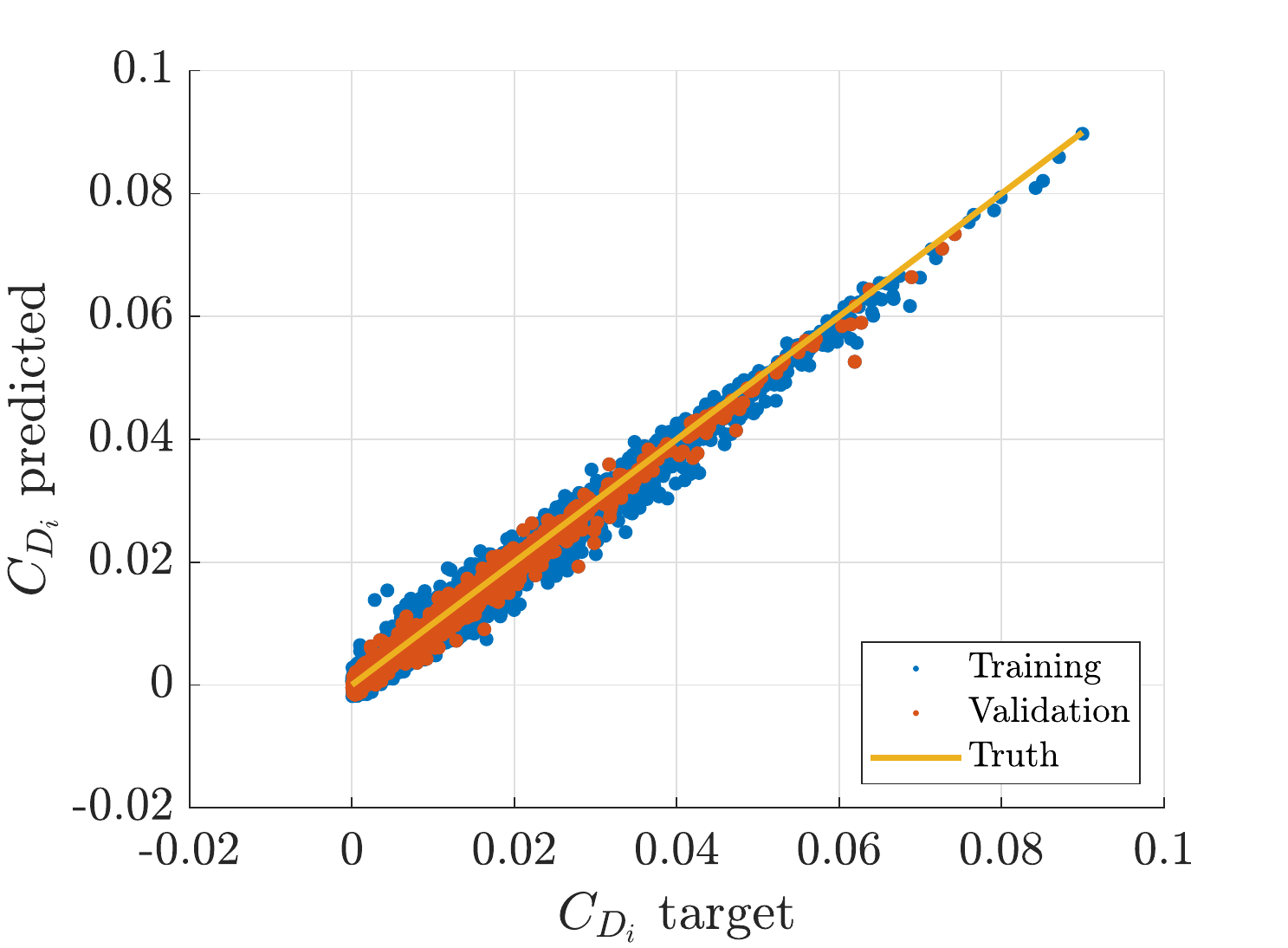}}
\caption{CNN predictions for $C_L$ and $C_{D_i}$}
\label{fig:train_val_pred}
\end{figure}

The CNN predictions for the training and validation sets are presented in Figs.~\ref{fig:train_val_pred} (a) and (b), for $C_L$ and $C_{D_i}$, respectively. For the lift coefficient, the training and validation predictions had $R^2=0.99$ and $R^2=0.98$, respectively. For the induced drag coefficients, the CNN predictions for the training and validation sets were $R^2=0.984$ and $R^2=0.975$, respectively.

Several of the hyperparameters of the CNN model were tuned to obtain satisfactory training convergence and accurate generalization over the validation set, with the fastest possible training time. These include the number of convolutional layers, the size of the convolutional kernels $\delta$, and the number of convolutional kernels $n_k$. These hyperparameters mainly affect the CNN model size and generalizable capacity. The other hyperparameters include the initial training rate, decay of the training rate, and total training iterations. These hyperparameters affect the training convergence. We only show the sensitivity of the CNN models to a few of the key hyperparameters here for brevity.

\begin{figure}[hbt!]
\centering
\subfloat[$C_L$]{\includegraphics[width=0.48\linewidth]{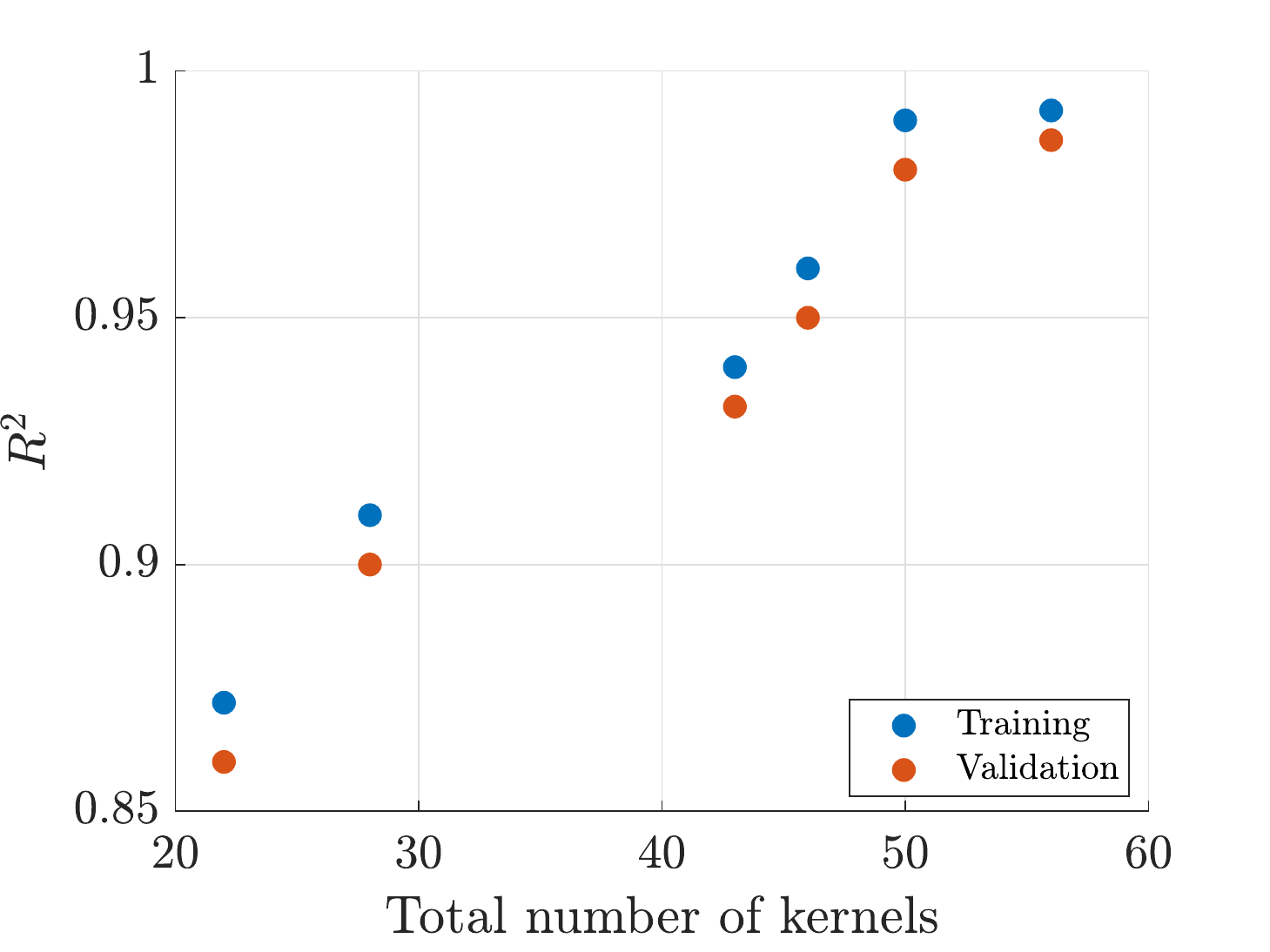}}
\subfloat[$C_{D_i}$]{\includegraphics[width=0.48\linewidth]{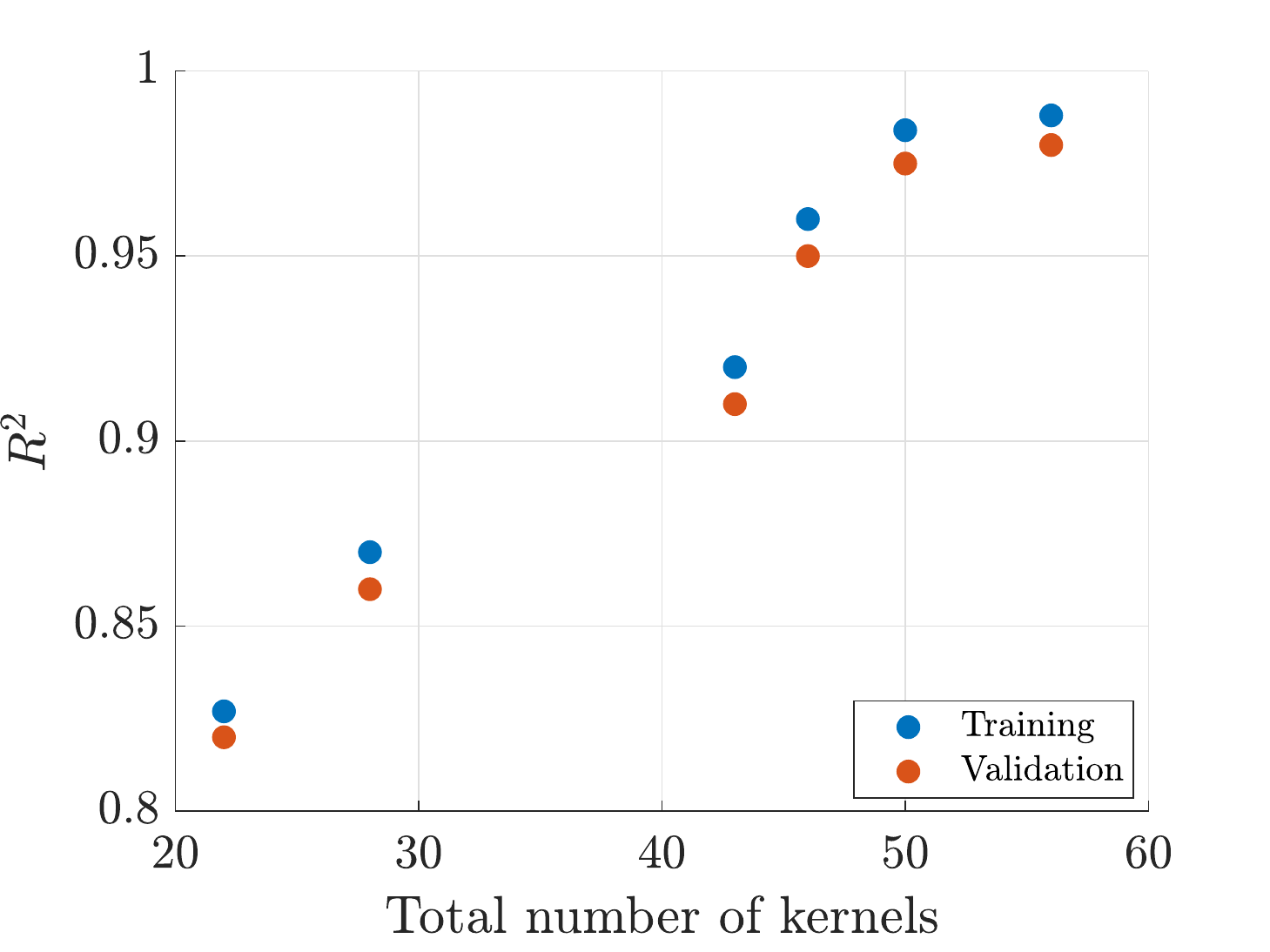}}\\
\caption{Sensitivity of training and validation accuracy to total number of kernels}
\label{fig:hyperparam_sens_ksize}
\end{figure}

First, we discuss the effect of the CNN model size (or the number of trainable parameters) on its predictive capability. We fix the number of convolutional layers as 4 and the kernel size $\delta$ in each layer. Thus, the CNN model size depends on the total number of kernels. The total number of kernels in the deep CNN model is varied by changing the number of kernels in each layer, $n_k$. It can be observed from Fig.~\ref{fig:hyperparam_sens_ksize} that on increasing the total number of kernels, the CNN's training and validation prediction accuracy kept improving, for both the $C_L$ and $C_D$, as indicated by their $R^2$ values. Validation and training prediction accuracy increased proportionally indicating that satisfactory training on 15,000 training sets usually assures reasonably generalizable prediction capability for other designs in the same parameter space. The network with 50 total kernels was selected as the optimal one. Further increasing the number of parameters to 56 resulted in a 30\% higher training cost but a very small improvement in the prediction accuracy. The $n_k$ for each layer of the various CNN models presented in Fig.~\ref{fig:hyperparam_sens_ksize}, is provided in Table \ref{tab:ksize_distribtion}. The model size for each of these CNN models is indicated by the total trainable parameters. The model size increases proportionally with the total number of kernels.

\begin{table}
\centering
\caption{$n_k$ distribution in each layer of the CNN models}
\label{tab:ksize_distribtion}       
\begin{tabular}{ccc}
\hline\noalign{\smallskip}
Total number of kernels & Kernel distribution & Total trainable parameters \\
\noalign{\smallskip}\hline\noalign{\smallskip}
22 & 4, 5, 6, 7 & 7829\\
28 & 4, 6, 8, 10 & 21365\\
43 & 5, 9, 13, 16 & 86973\\
46 & 5, 10, 14, 17 & 110646\\
50 & 5, 10, 15, 20 & 139451\\
56 & 6, 11, 16, 21 & 206164\\
\noalign{\smallskip}\hline
\end{tabular}
\end{table}

\begin{figure}[hbt!]
\centering
\subfloat[$C_L$]{\includegraphics[width=0.48\linewidth]{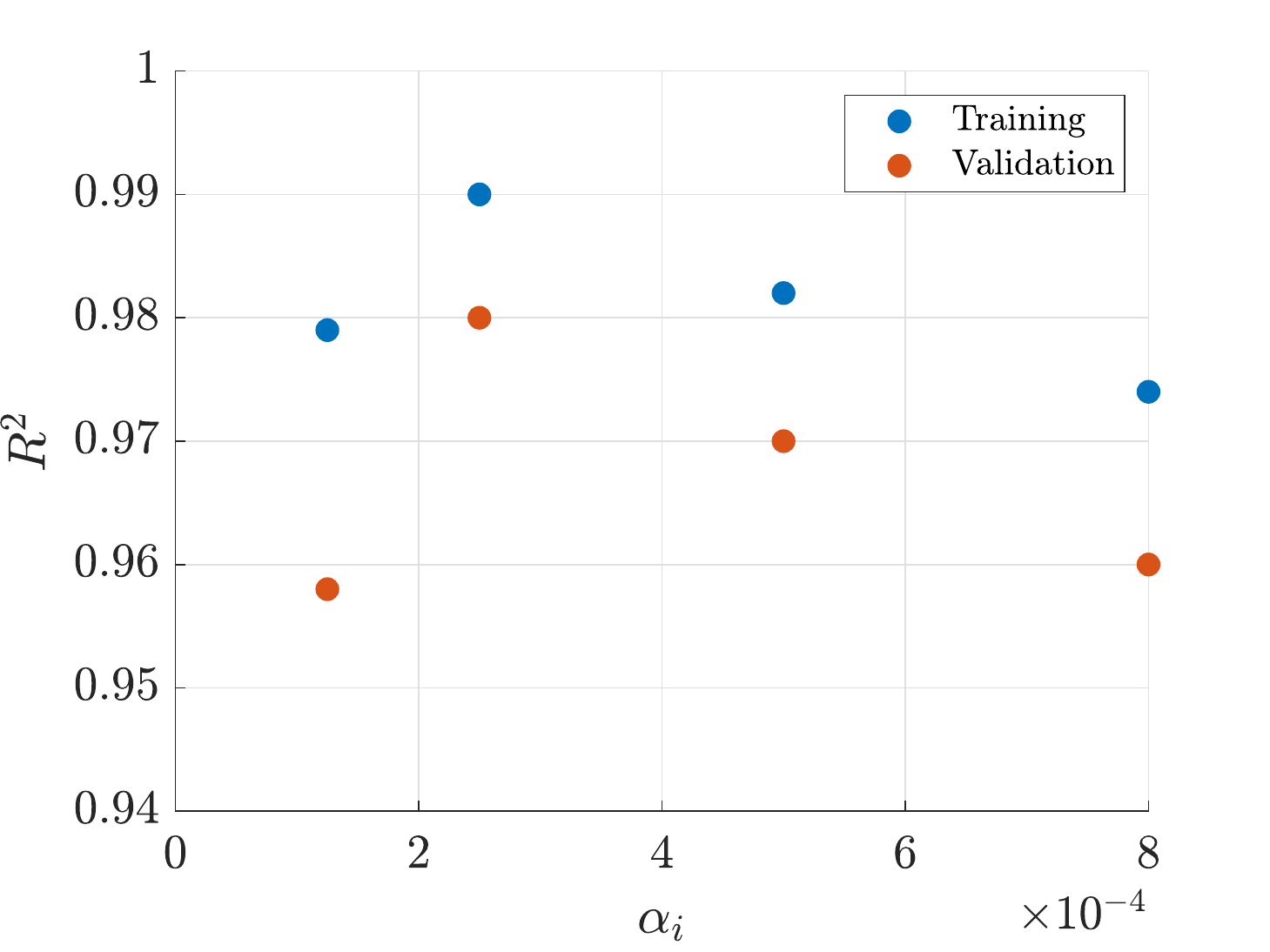}}
\subfloat[$C_{D_i}$]{\includegraphics[width=0.48\linewidth]{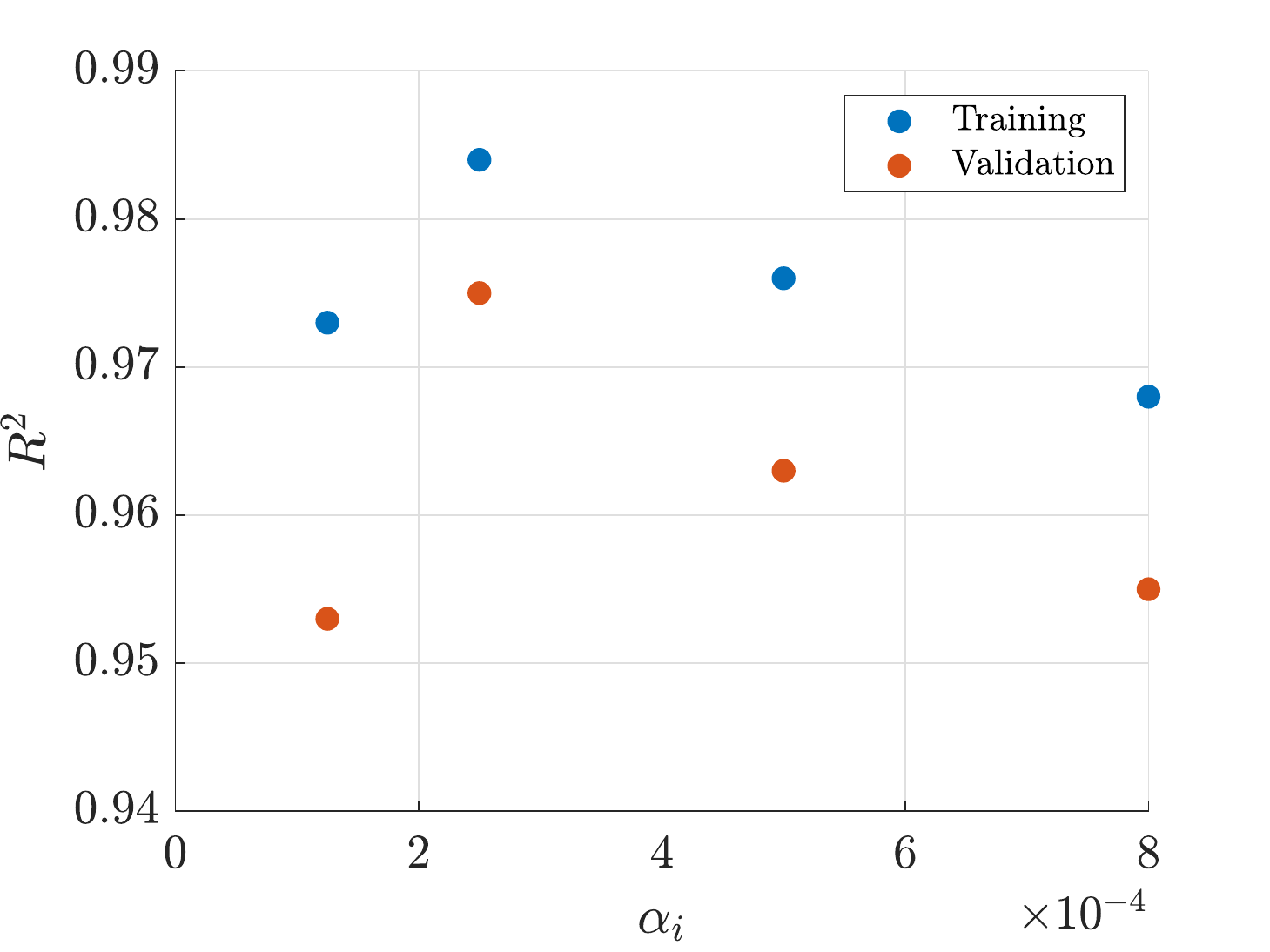}}\\
\caption{Sensitivity of training and validation accuracy to $\alpha_i$}
\label{fig:hyperparam_sens_alpha}
\end{figure}

The initial learning and the learning rate decay were also tuned to enable maximum training and validation accuracy with minimum possible training cost. Initial studies suggested 60,000 iterations to suffice. Thus, we select several initial learning rates $\alpha_i$ and adjusted our linear learning rate decay rate per iteration to reach the same final learning rate after 60,000 iterations. The results presented in Fig.~\ref{fig:hyperparam_sens_alpha} indicate that starting with $\alpha_i=0.00025$ is optimal. Training with a higher initial rate leads to higher oscillations in the stochastic gradient descent and lower convergence. On the other hand, a slower initial learning rate leads to a slow convergence rate by the time 60,000 iterations have been completed. The $C_L$ and $C_{D_i}$ predicting networks with optimal hyperparameters were used subsequently for the optimization. The complete list of hyperparameters of the optimal network is shown below in Table \ref{tab:hyperparam_tunedCNN}. The kernel size and the number of kernels in each layer of the CNN model are shown.

\begin{table}
\centering
\caption{Hyperparameters of the optimally tuned CNN model}
\label{tab:hyperparam_tunedCNN}       
\begin{tabular}{lc}
\hline\noalign{\smallskip}
Hyperparameters & Values \\
\noalign{\smallskip}\hline\noalign{\smallskip}
Kernel size ($\delta$) & 4, 3, 3, 3\\
Number of kernels ($n_k$) & 5, 10, 15, 20\\
Initial learning rate ($\alpha_i$) & 2.5e-4 \\
Learn rate decay factor & 0.999 \\
Learn rate decay period & 1 epoch \\
Epochs & 3000 \\
Batch size & 750 \\
\noalign{\smallskip}\hline
\end{tabular}
\end{table}

\begin{figure}[hbt!]
\centering
\subfloat[$C_L$]{\includegraphics[width=0.49\linewidth]{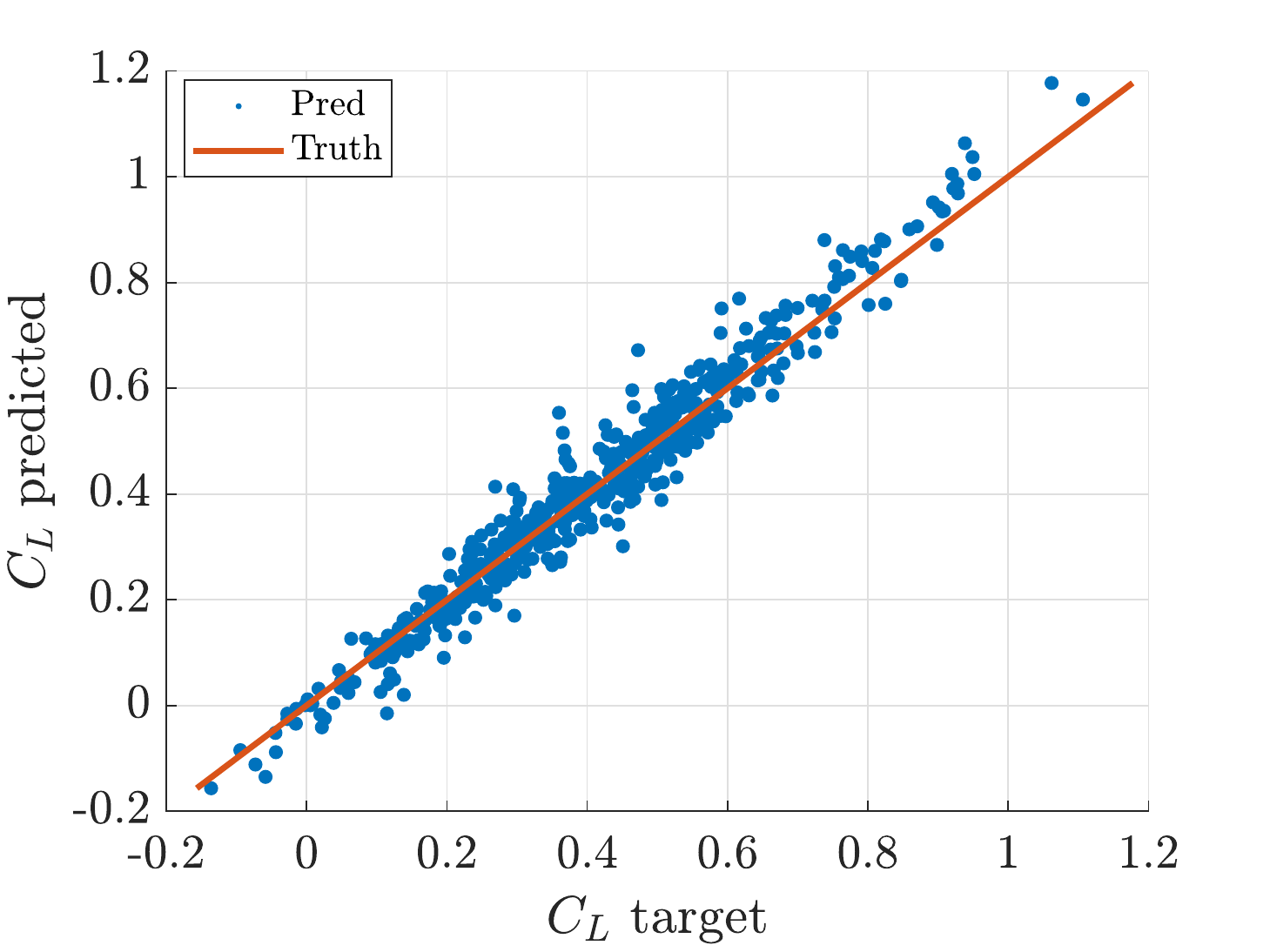}}
\subfloat[$C_{D_i}$]{\includegraphics[width=0.49\linewidth]{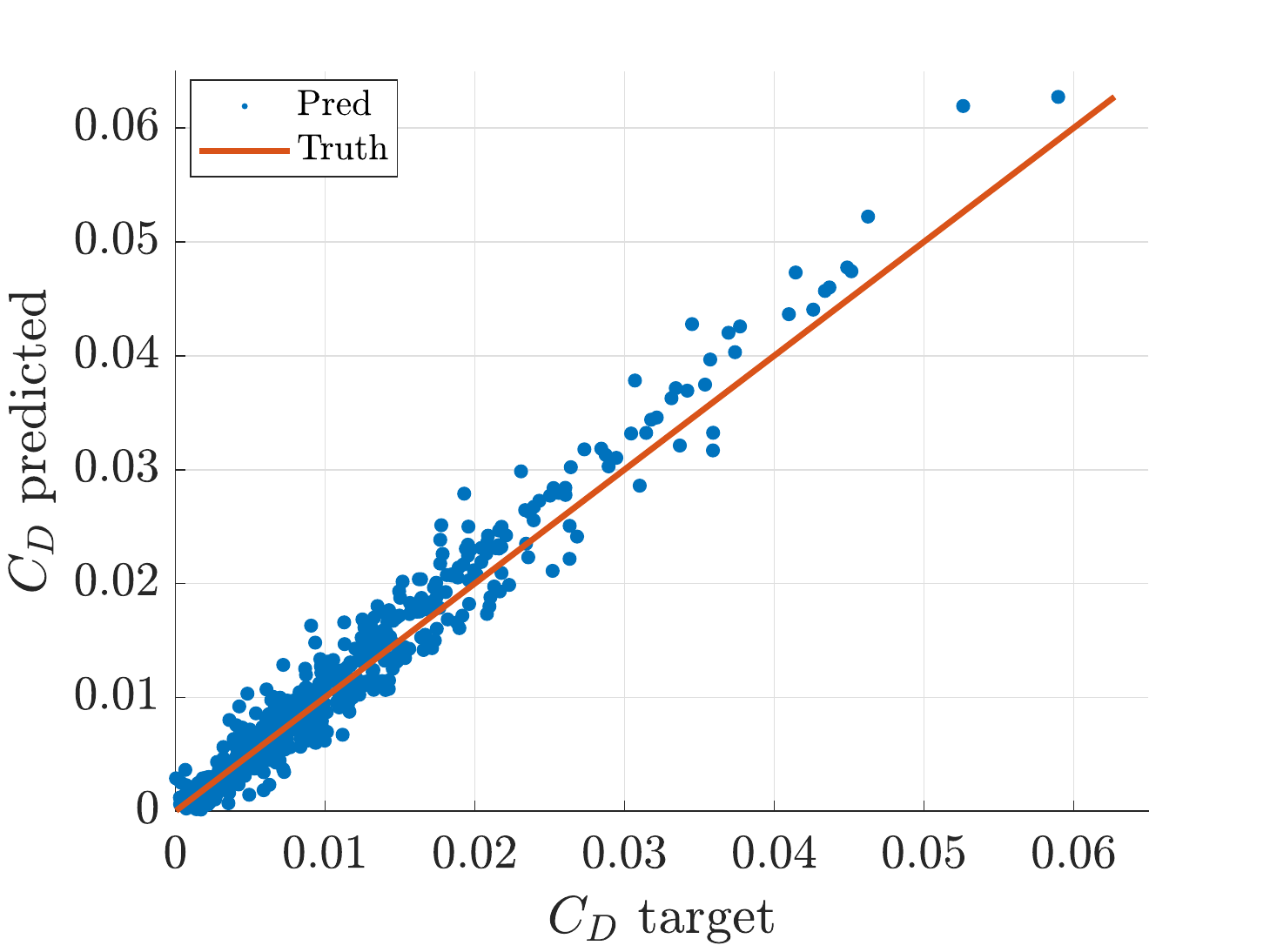}}
\caption{CNN test predictions for $C_L$ and $C_{D_i}$}
\label{fig:test_pred}
\end{figure}

We finally test our CNN on the test data set. The test data set will also serve as the initial starting designs for the subsequent gradient-based optimization study. The results presented in Fig.~\ref{fig:test_pred} (a) and (b) show the test predictions of the optimally tuned CNN network for the $C_L$ and $C_{D_i}$, respectively. The lift coefficient predictions show reasonable agreement with the ground truth results. The drag predictions show excellent agreement at low drag values and deviations from the true values are only observed for a few outliers at high values of the induced drag coefficient. Such high induced drag values will rarely be encountered when we optimize our designs to reduce induced drag and thus can be safely neglected for the present study. The predicted $C_L$ and $C_{D_i}$ have $R^2=0.945$ and $R^2=0.95$, respectively. 

\subsection{Deep CNN feature space exploration}
\begin{figure}[hbt!]
\centering
\subfloat[wing geometry]{\includegraphics[width=0.52\linewidth]{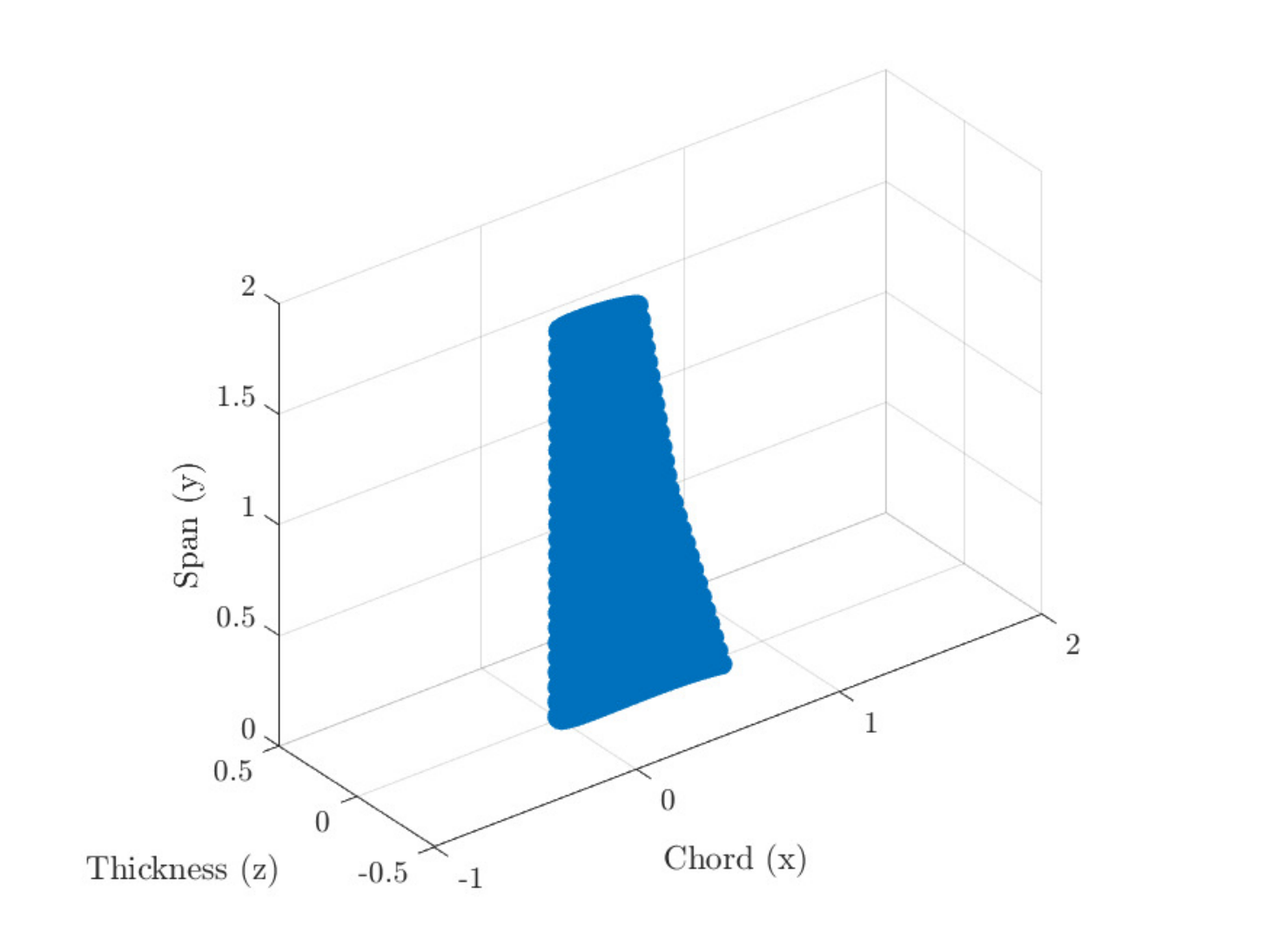}}
\subfloat[wing level-set]{\includegraphics[width=0.52\linewidth]{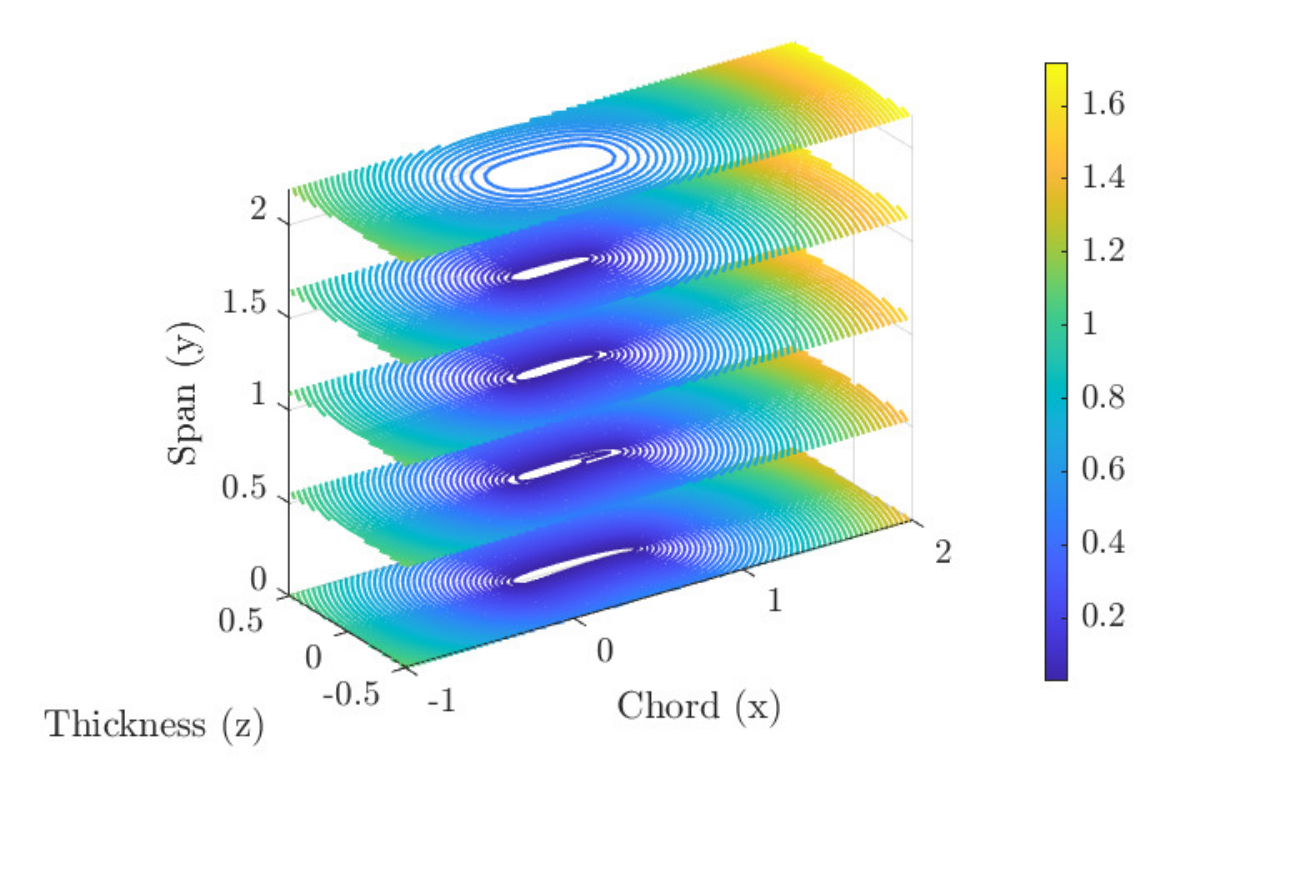}}
\caption{Representative wing geometry and level-set}
\label{fig:representative_wing}
\end{figure}

We have shown earlier that the CNN can learn a functional relation between the shapes and their aerodynamic force coefficients. Here we will explore the feature maps of various convolutional layer kernels to interpret the learning mechanism facilitated by the hierarchical structure of the deep CNN. We first present a representative wing and its corresponding level-sets in Figs.~\ref{fig:representative_wing} (a) and (b), respectively. These will be used for generating the deep CNN feature maps.

\begin{figure}[hbt!]
\centering
\subfloat[Kernel 1 features]{\includegraphics[width=0.48\linewidth]{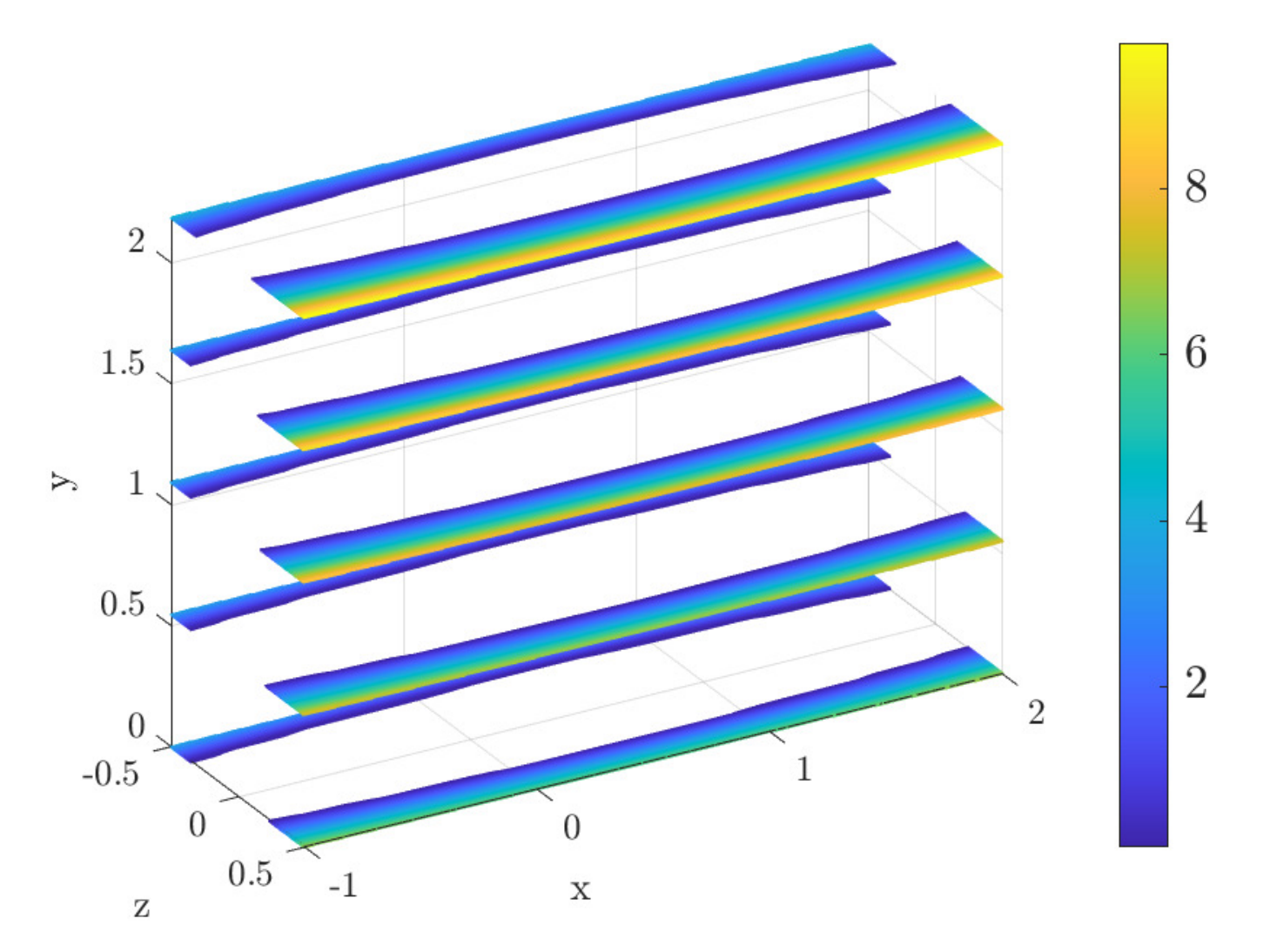}}
\subfloat[Kernel 2 features]{\includegraphics[width=0.48\linewidth]{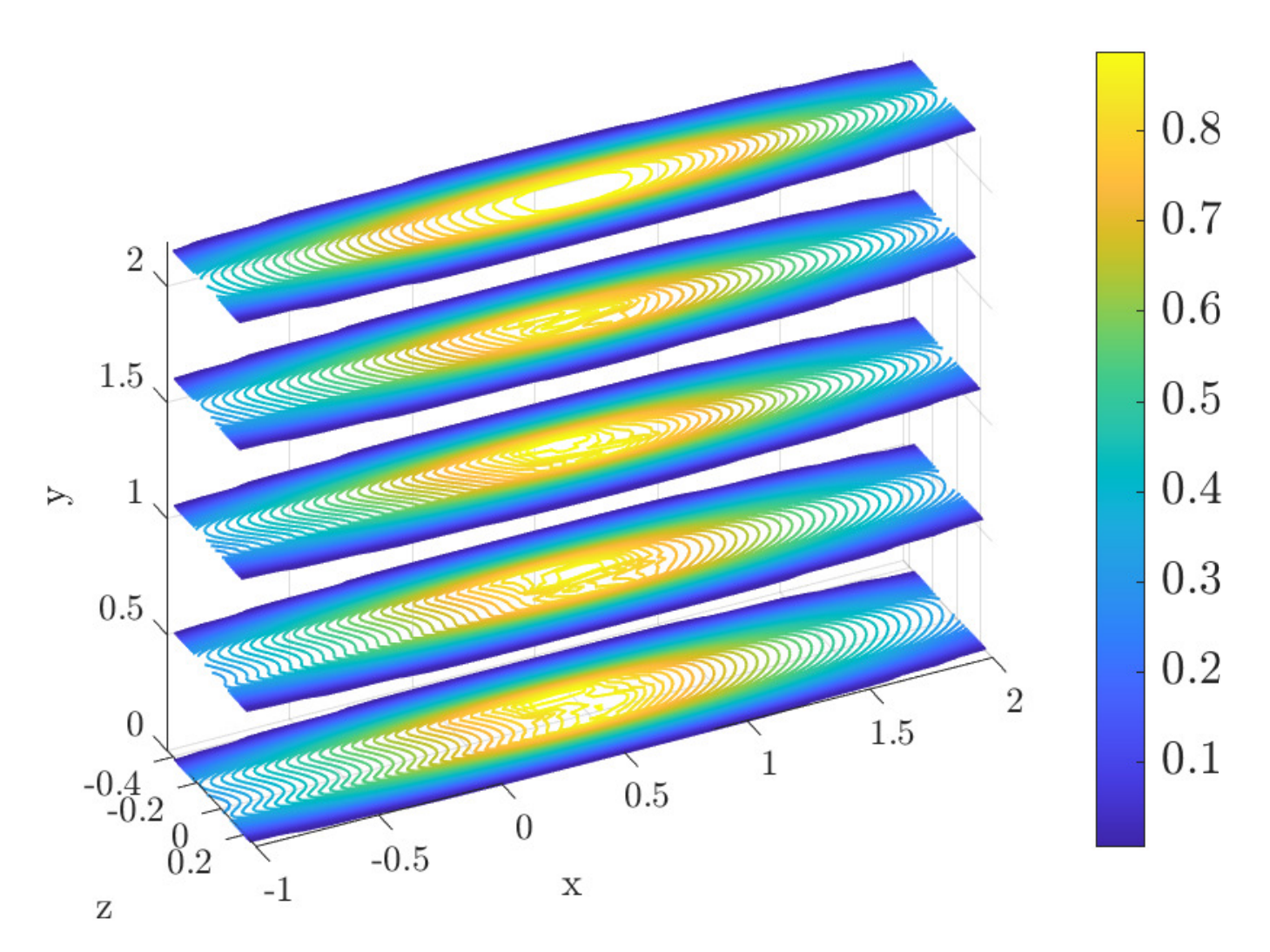}}\\
\subfloat[Kernel 3 features]{\includegraphics[width=0.48\linewidth]{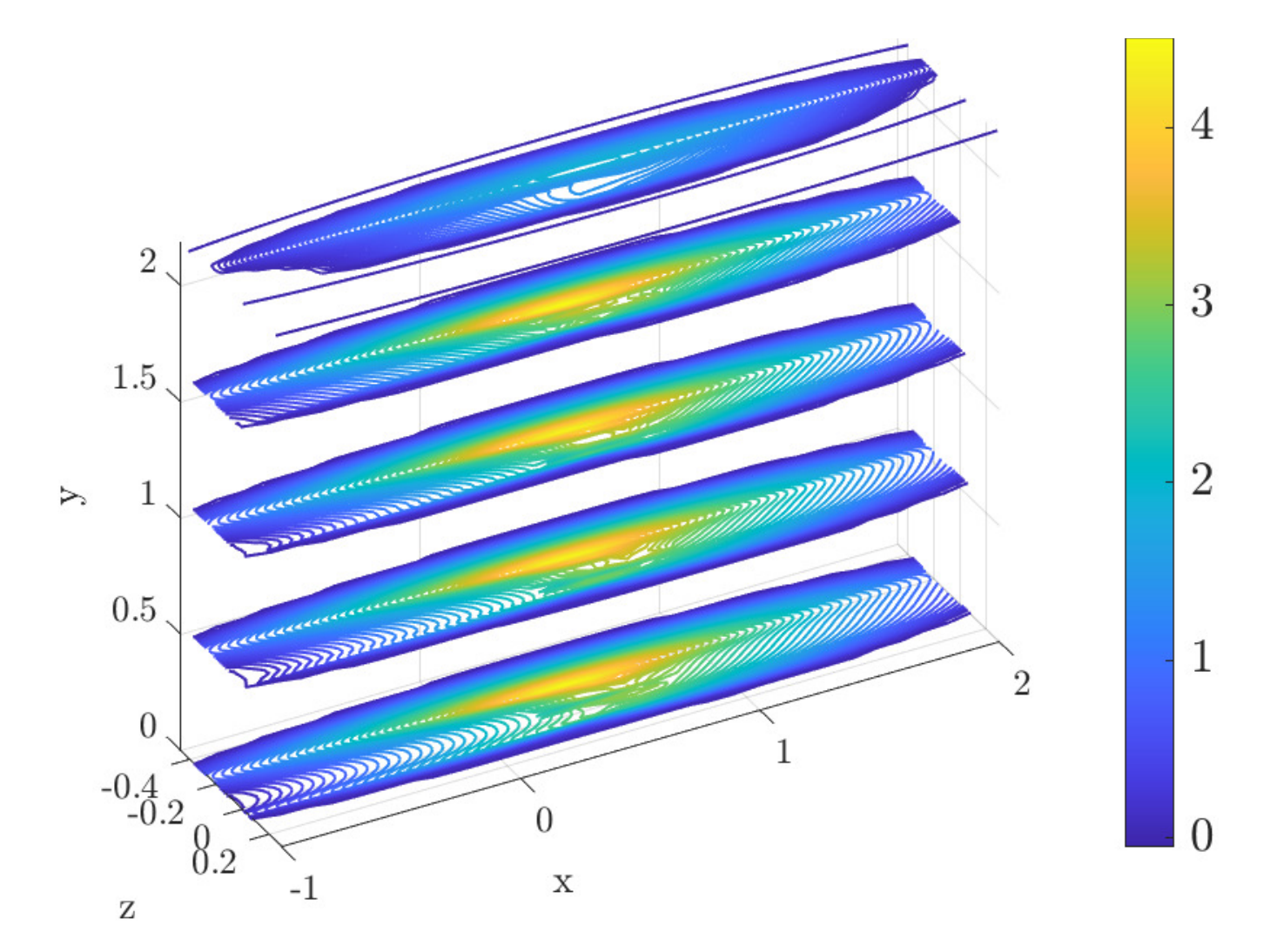}}
\caption{Three-dimensional feature maps: first convolutional layer kernels}
\label{fig:CNNlayer1out}
\end{figure}

We first explore the feature maps of the first convolutional layer of the network. Figs.~\ref{fig:CNNlayer1out} (a), (b) and (c) present the feature maps of the first, second and third kernels, respectively. Various slices of the feature maps are presented along the $y$ axis of feature maps to show how they vary along the $y$-axis. Each slice shows the behavior of the feature maps along the $x-z$ plane at a specific $y$ value. Fig.~\ref{fig:CNNlayer1out} (a) indicates that the first kernel focuses on the feature maps of the wing level-set (Fig.~\ref{fig:representative_wing} (b)) lying near the upper limit of the $z$ axis. However, the weights for these kernel feature maps increase as move upward along the $y$ axis. The second kernel (Fig.~\ref{fig:CNNlayer1out} (b)) focuses on the input features lying near the center of the $x-z$ plane. The feature maps are almost constant along the $y$-axis. The third kernel (Fig.~\ref{fig:CNNlayer1out} (c)) also focuses on the features near the center of the $x-z$ plane but more towards the lower limit of the $y$ axis (the root of the wing). The kernel feature maps presented here are analogous to three-dimensional modes learned by the first convolutional layer as a result of the decomposition of input level-set features into these kernels. However, unlike linear projection methods like POD, these CNN kernel feature maps are not orthogonal in nature and decomposed via a nonlinear process.

\begin{figure}[hbt!]
\centering
\subfloat[Kernel 2 features]{\includegraphics[width=0.48\linewidth]{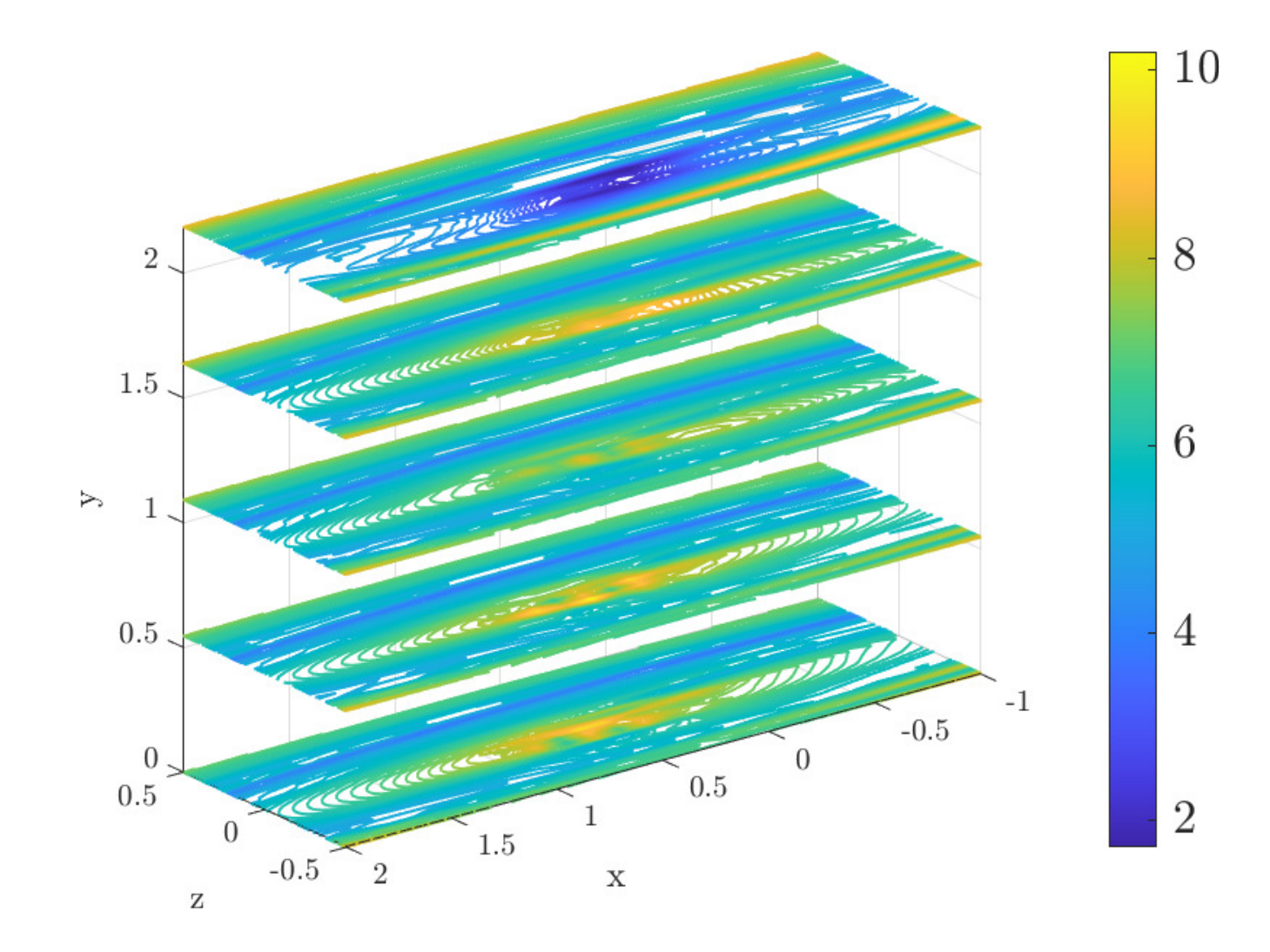}}
\subfloat[Kernel 7 features]{\includegraphics[width=0.48\linewidth]{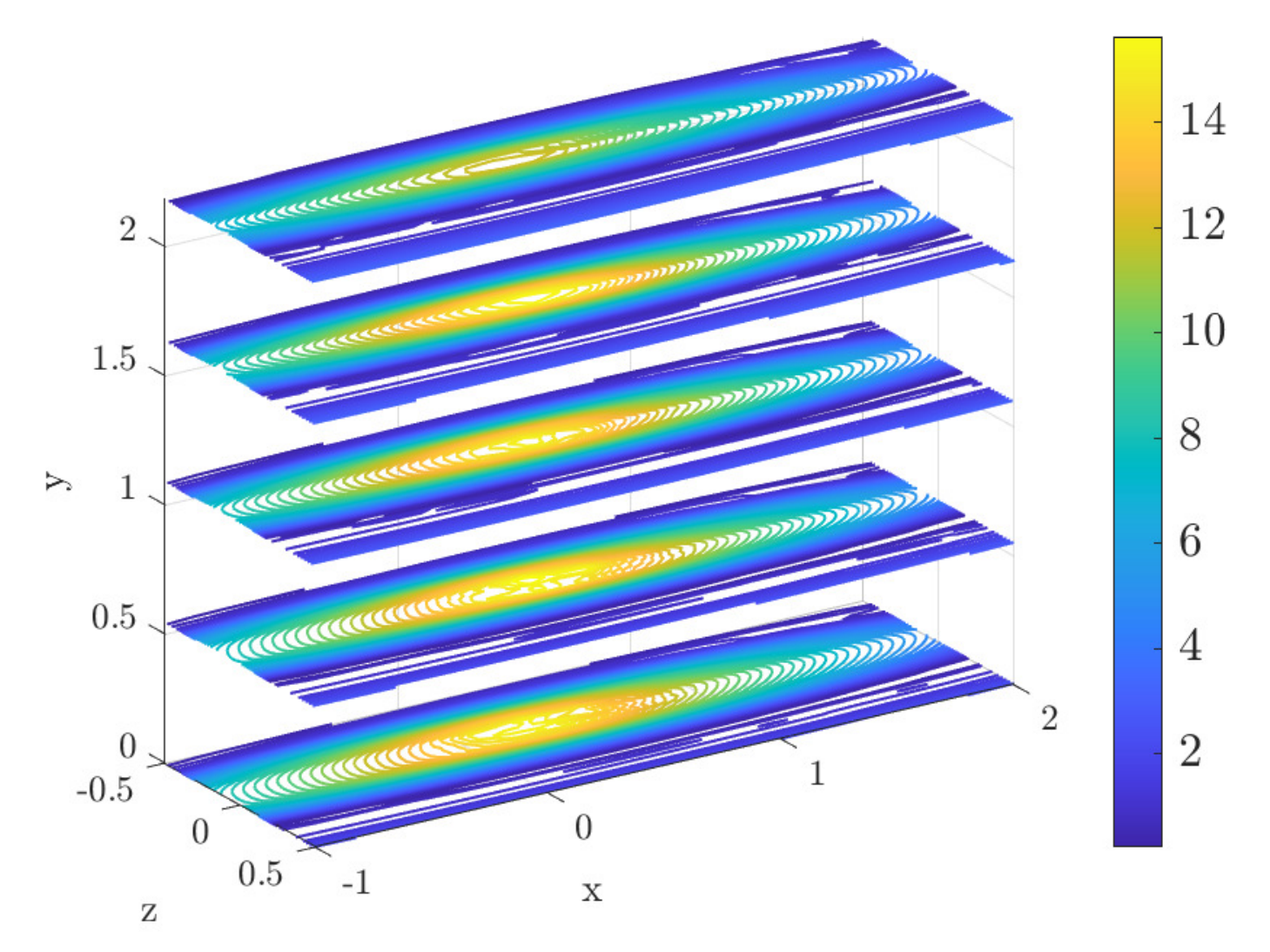}}\\
\subfloat[Kernel 14 features]{\includegraphics[width=0.48\linewidth]{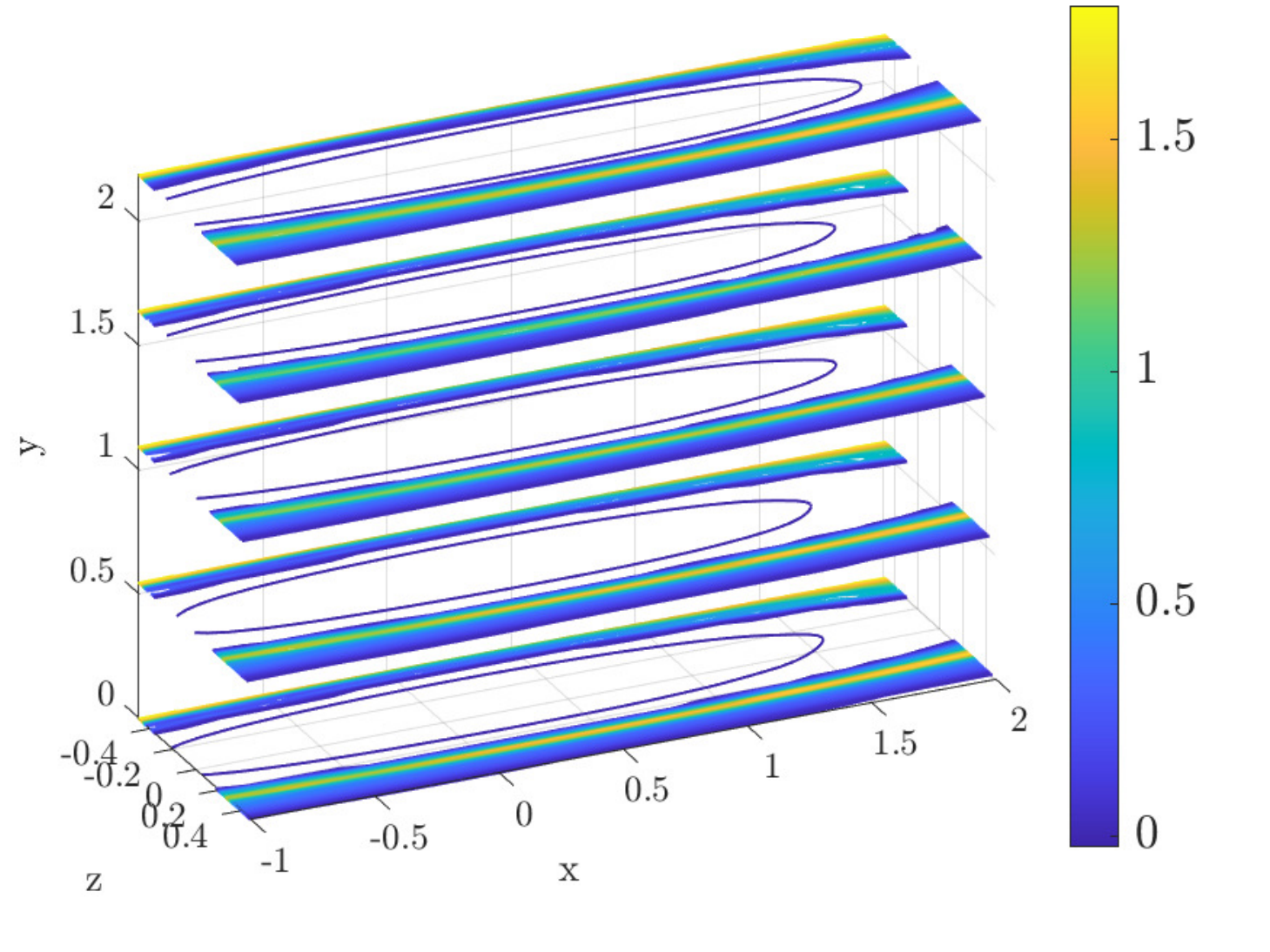}}
\caption{Three-dimensional feature maps: third convolutional layer kernels}
\label{fig:CNNlayer3out}
\end{figure}

Similar to the first convolutional layer, we also investigate the features maps of various kernel outputs of the third convolutional layer. We specifically present the features of the 2nd, 7th and 14th kernel outputs of the third convolutional layer in Figs.~\ref{fig:CNNlayer3out} (a), (b) and (c), respectively. The 2nd kernel output shows a feature map whose $x-z$ features vary non-uniformly along the $y$-axis. The feature map (Fig.~\ref{fig:CNNlayer3out} (a)) focuses on the features near the center of the $x$-$z$ plane for lower values of $y$. It shows a uniform behavior along the $x$-$z$ plane near the middle of the $y$-axis and again shows higher weights near the center of the $x$-$z$ plane with a further increase in $y$. Finally, near the upper limit of the $y$ axis, the feature map weights along the $x-z$ plane progressively increase away from the center. For the 7th kernel, the feature maps (Fig.~\ref{fig:CNNlayer3out} (b)) focus specifically on the features close to the center of the $x-z$ and such output features show the highest weights near the lower limit of the $y$ axis. The $x-z$ features show lower concentration near the center as we move up along the $y$ axis. The 14th kernel output features (Fig.~\ref{fig:CNNlayer3out} (c)) focus on a specific set of level-set isocurves lying midway between the center and the outer limit of the $x$-$z$ plane. The features remain almost constant along the $y$ axis.

Overall, we observe a hierarchy in the feature maps of the kernels depending on which layer they are located in. The feature maps of kernels in the outer layers learn global features. Whereas, the kernel feature maps in the deeper layers capture the more localized features showing widely varying distribution along the various axes. This hierarchy follows the hierarchical structure of the deep CNN as there are 15 kernels in the third, 10 kernels in the second layer, and 5 in the first layer. Thus, the nonlinear decomposition of the 10 kernel outputs of the second convolutional layer outputs into the 15 kernel features of the third layer is a far more complex phenomenon compared to a decomposition of the input level-sets to the 5 kernel feature maps of the first convolutional layer. The deep CNN model is a compositional representation of all the kernels. Thus, the widely varying kernel feature maps capturing shape information from global to local scales, provide stability to the composite deep CNN map against shape deformations occurring at various scales. Such geometric assumptions of the deep CNN framework were abstractly represented earlier. Here we provide physical interpretation by exploring the kernel feature maps of the various convolutional layers.

\subsection{Gradient-based optimization results}
\begin{figure}[ht!]
\centering
\includegraphics[width=0.98\linewidth]{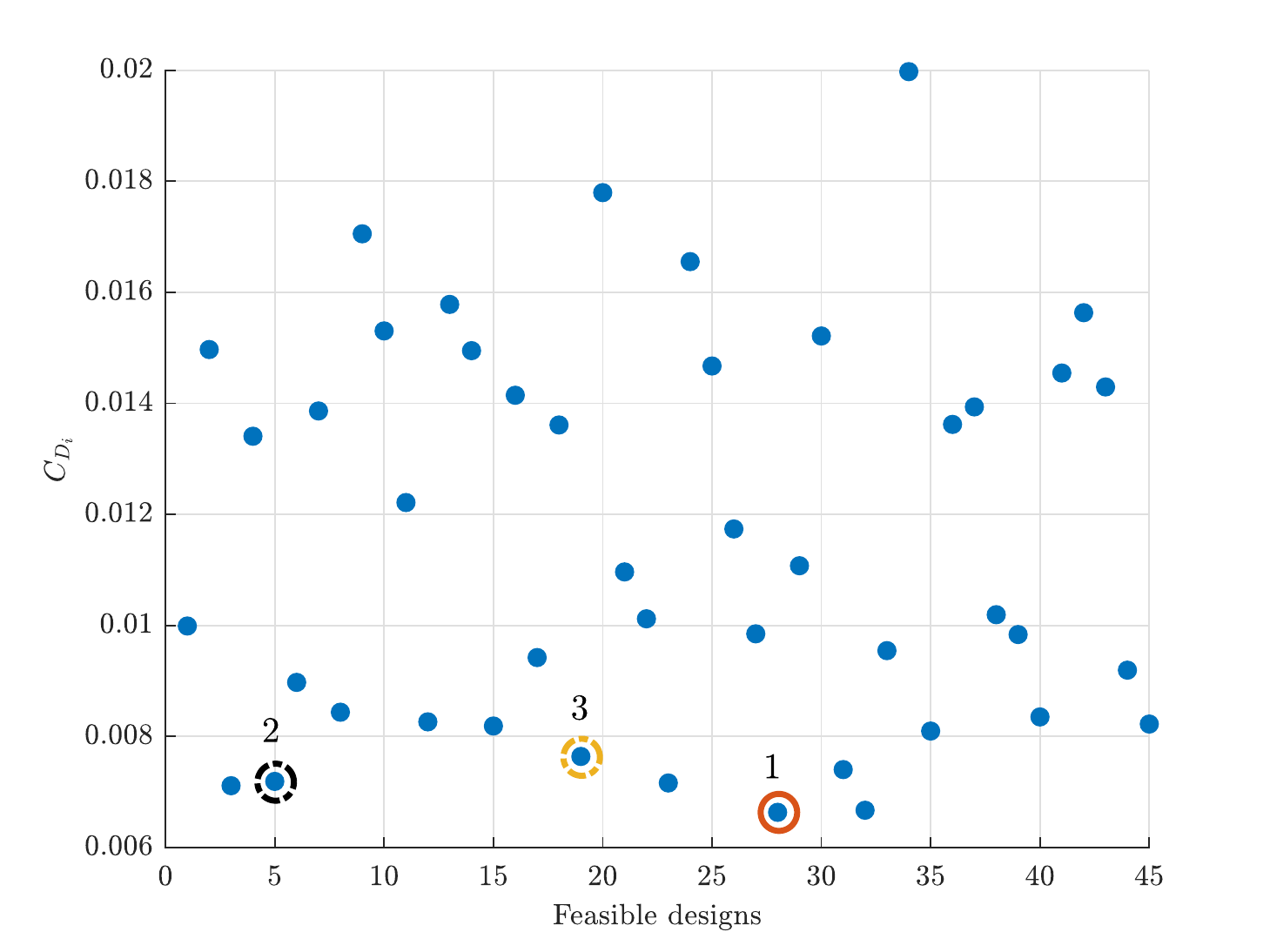}
\caption{Feasible designs obtained via DL-ROM-based optimization}
\label{fig:feasible_des}
\end{figure}

We perform our gradient-based optimization by starting with several different initial design points. Such design space exploration is performed by assuming the possibility for the objective function to show non-convex behavior in the thirty-three dimensional design space. Hence 50 initial designs are selected randomly from the test set presented earlier and the DL-ROM-based optimization was performed for each of them. The constrained gradient-based optimization problem is solved via Sequential Quadratic Programming (SQP) algorithm from Matlab's \textit{fmincon} function \cite{matlabopt}. The optimization could not converge to a feasible design for certain initial designs. However, for all the others a feasible design could be obtained at the end of the optimization process. These are presented in Fig.~\ref{fig:feasible_des}. As previously suspected, the final objective values of the feasible designs vary widely indicating the non-convex nature of the objective function. For the feasible designs shown in Fig.~\ref{fig:feasible_des}, the objective functions are obtained with $R^2=0.94$ compared to their counterparts computed with the FOM. This demonstrates the accuracy of the CNN-based DL-ROM during the optimization process.

\begin{figure}[ht!]
\centering
\includegraphics[width=0.96\linewidth]{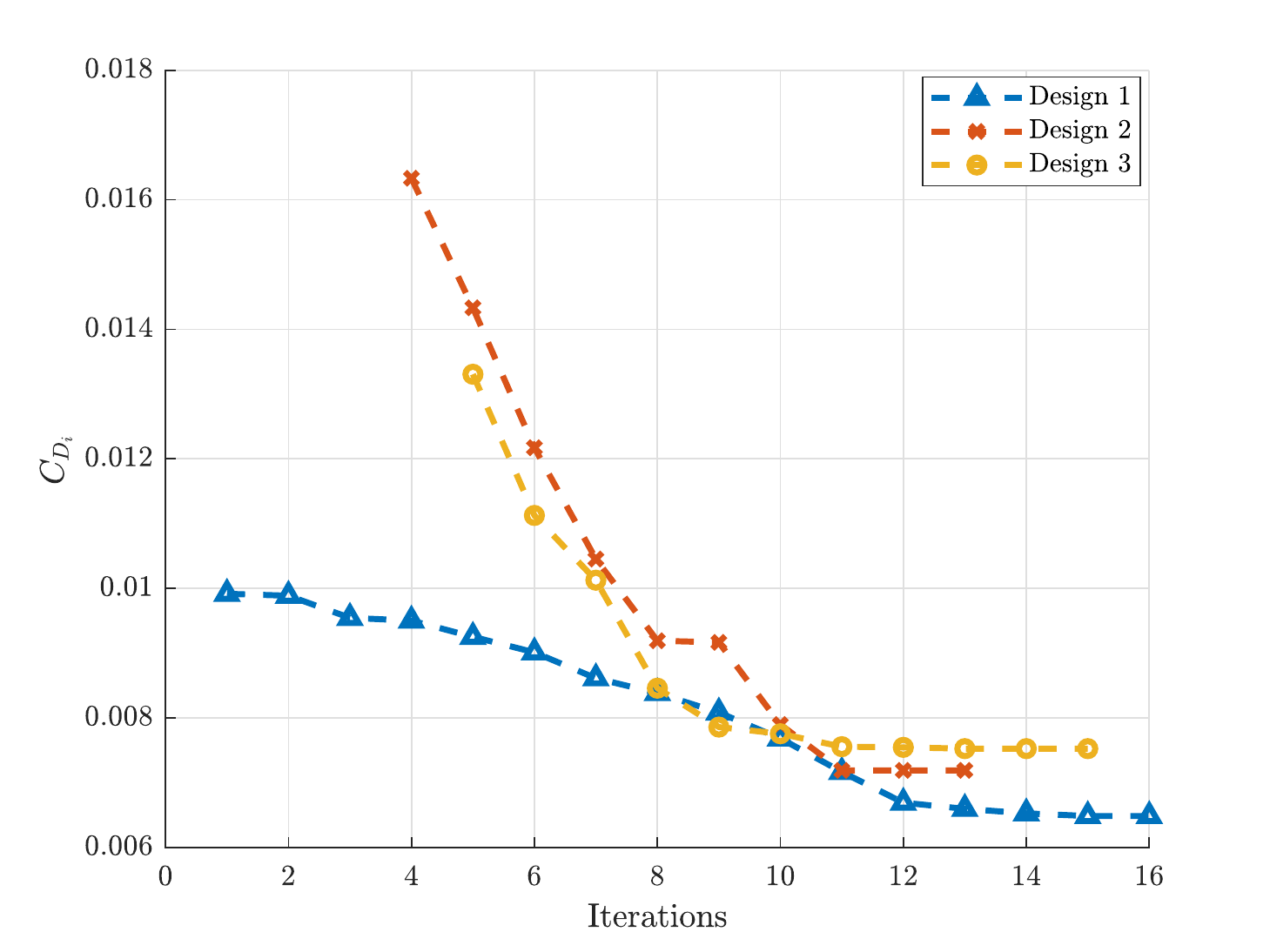}
\caption{Convergence history of gradient-based optimization for a 3D wing problem}
\label{fig:SQP_conv_hist}
\end{figure}

Since we have employed a gradient-based optimization, each converged feasible design in Fig.~\ref{fig:feasible_des} can be considered as a local optimum. From these local optima we select three designs, 1, 2 and 3 (marked with circles in Fig.~\ref{fig:feasible_des}), for further inspection. Design 1 has an objective of $C_{D_i}=0.0066$, which is the minimum amongst all the feasible designs obtained. Therefore, it can be considered the best design obtained from the optimization. Design 2 has $C_{D_i}=0.0072$ and design 3 has $C_{D_i}=0.0076$. These are considered as two representative feasible designs obtained from the optimization process, having an objective value close to the best design obtained. The gradient-based optimization history for designs 1, 2 and 3 are presented in Fig.~\ref{fig:SQP_conv_hist}. Intermediate objective values are only shown for those iterations where the designs have satisfied all the constraints. We can see that for design 1, only the initial design did not satisfy all the constraints. However, for design 3, the feasible designs could be obtained only after the fifth iteration. This shows that the DL-ROM-based constrained optimization process is capable of starting from arbitrary, infeasible designs and converging to feasible local minima for most cases.

\begin{figure}[hbt!]
\centering
\includegraphics[width=0.96\linewidth]{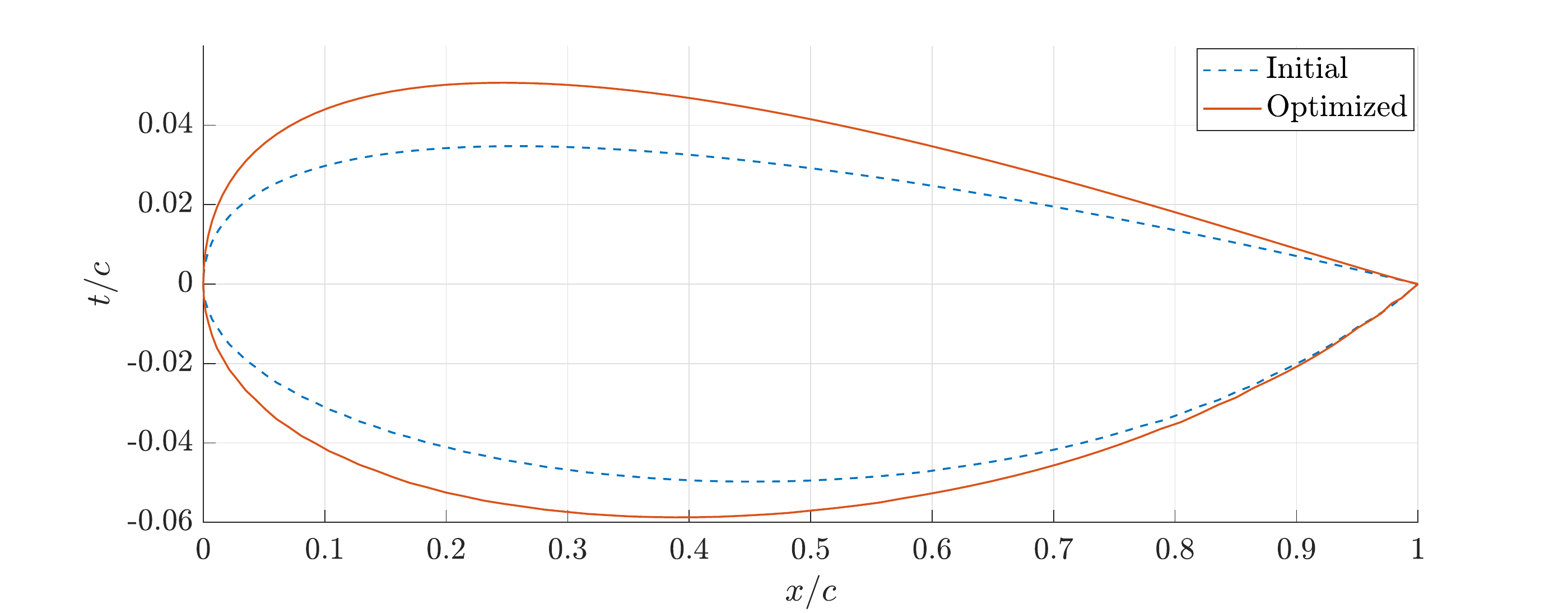}
\caption{Comparison of initial and optimized design 2 wing sections: root}
\label{fig:wing_sec_comp_root}
\end{figure}

\begin{figure}[hbt!]
\centering
\includegraphics[width=0.96\linewidth]{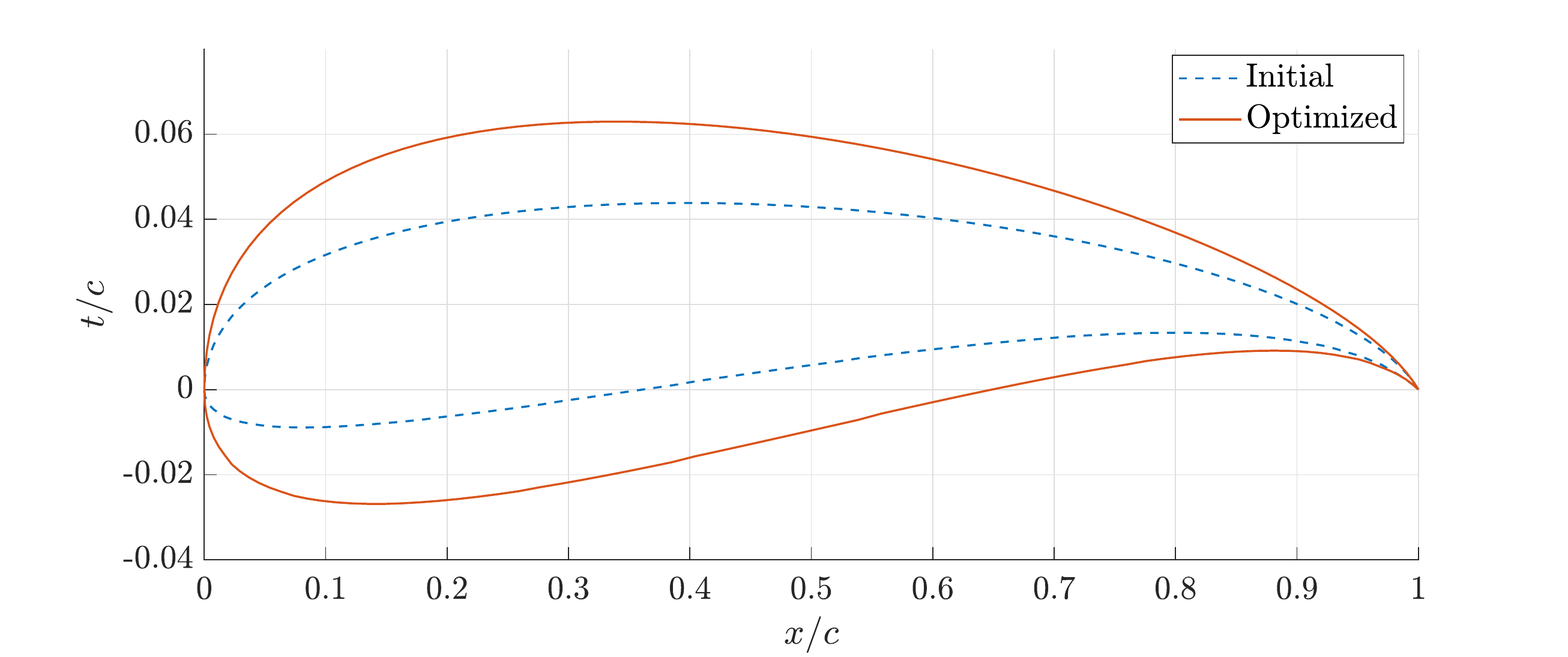}
\caption{Comparison of initial and optimized design 2 wing sections: 33\% span}
\label{fig:wing_sec_comp_p33span}
\end{figure}

\begin{figure}[hbt!]
\centering
\includegraphics[width=0.96\linewidth]{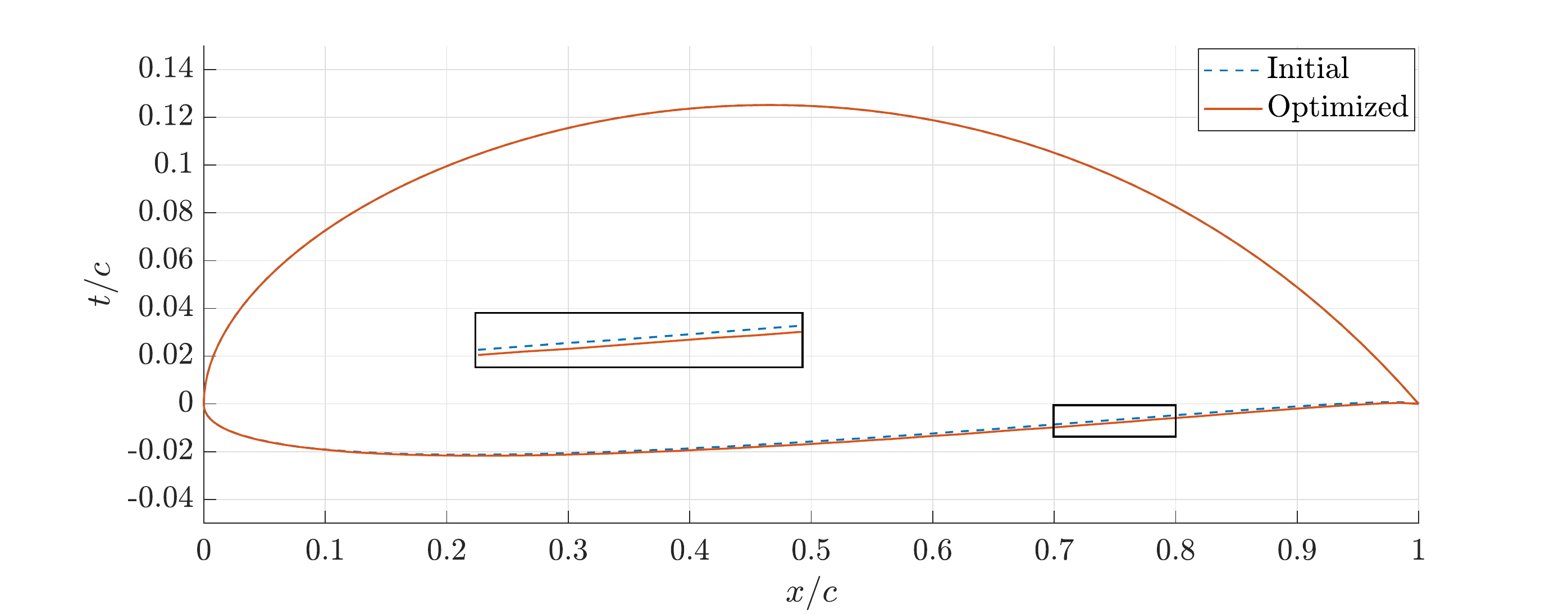}
\caption{Comparison of initial and optimized design 2 wing sections: 67\% span}
\label{fig:wing_sec_comp_p67span}
\end{figure}

\begin{figure}[hbt!]
\centering
\includegraphics[width=0.96\linewidth]{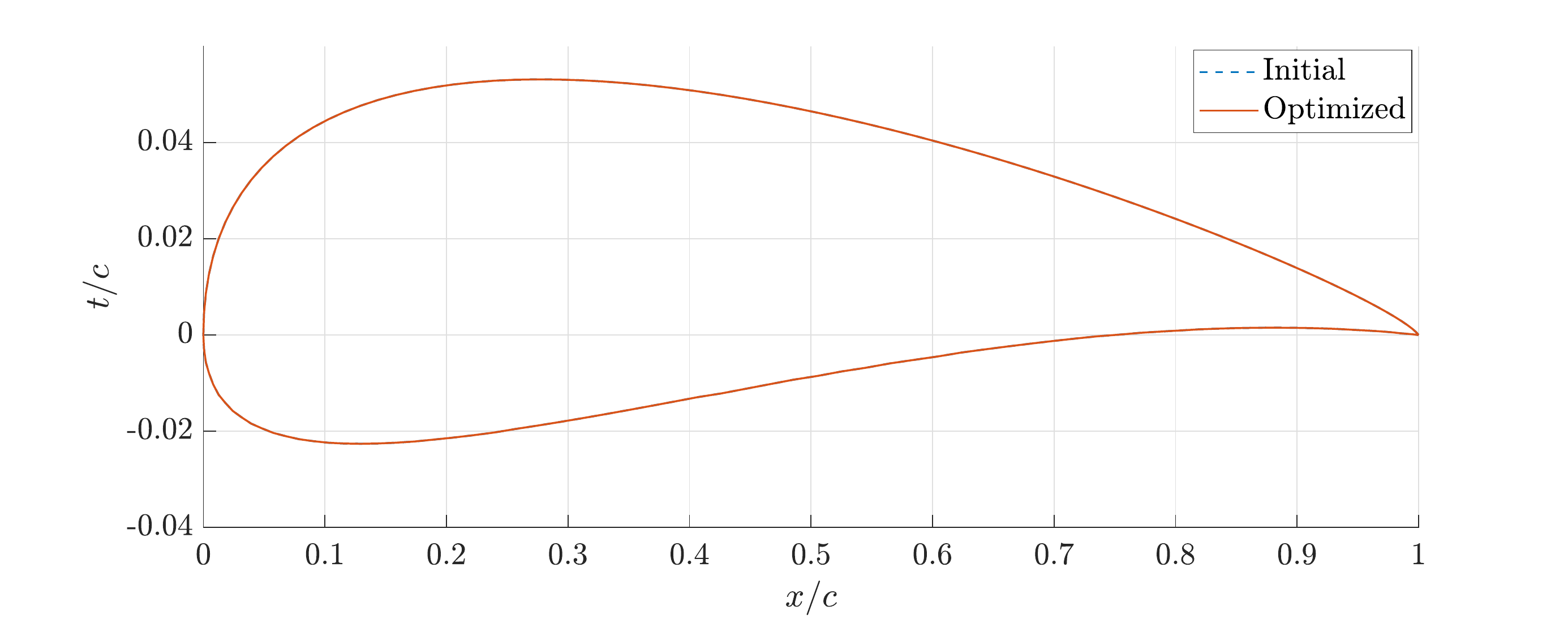}
\caption{Comparison of initial and optimized design 2 wing sections: tip}
\label{fig:wing_sec_comp_tip}
\end{figure}

To further illustrate the effects of the optimization, we investigate the differences between the initial and optimized wing configurations for one of the representative designs, design 2. For the comparison of the wing shapes, the wing sections for the initial and optimized configuration are obtained at four different locations along the span: the root, 33\%, and 67\% of the span, and the wingtip. These are shown in Figs.~\ref{fig:wing_sec_comp_root}, \ref{fig:wing_sec_comp_p33span}, \ref{fig:wing_sec_comp_p67span} and \ref{fig:wing_sec_comp_tip}, respectively. We can observe that the optimization process significantly changed the wing's thickness distribution along the span. The thickness at the wing root increased, especially on the upper surface (Fig.~\ref{fig:wing_sec_comp_root}). A significant increase in the wing section thickness is also observed at 33\% of the span (Fig.~\ref{fig:wing_sec_comp_p33span}). A small increase in the wing sectional thickness at 70-80\% of the chord length on the lower surface is observed (Fig.~\ref{fig:wing_sec_comp_p67span}) at 67\% of the span. The inset showing a magnification is provided to illustrate this. No notable changes were observed at the wing root section. Apart from these changes in the wing sectional thickness, the wing's incident angle of attack also changed from $5.8$ degrees to $4.9$ degrees due to the optimization process. These changes first ensured that the wing satisfies the thickness and the $C_L$ constraints as shown in Fig.~\ref{fig:SQP_conv_hist}, and then the induced drag reached a local minimum.

\begin{figure}[hbt!]
\centering
\subfloat[$C_L$ distribution: initial]{\includegraphics[width=0.49\linewidth]{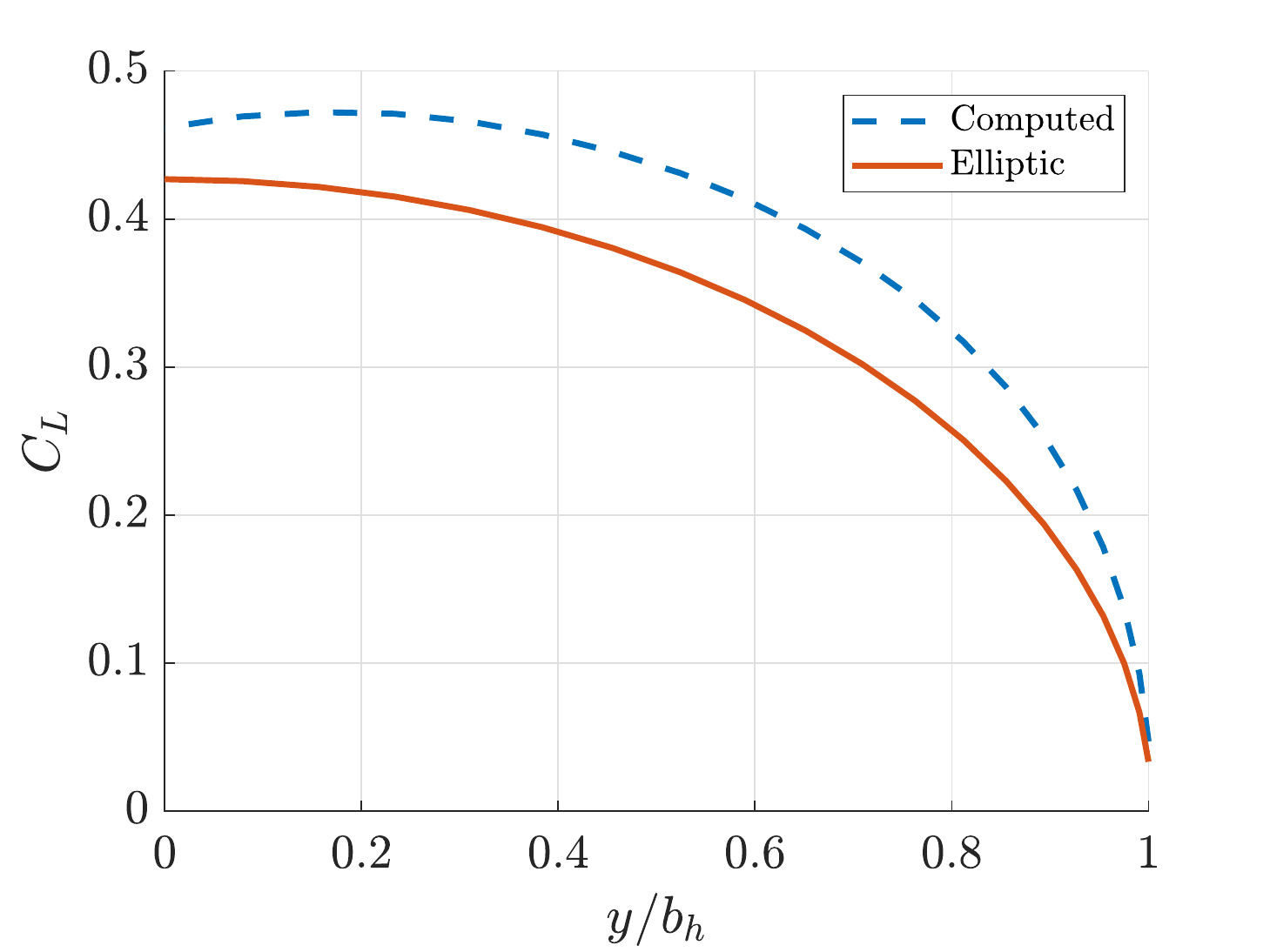}}
\subfloat[$C_L$ distribution: optimized]{\includegraphics[width=0.49\linewidth]{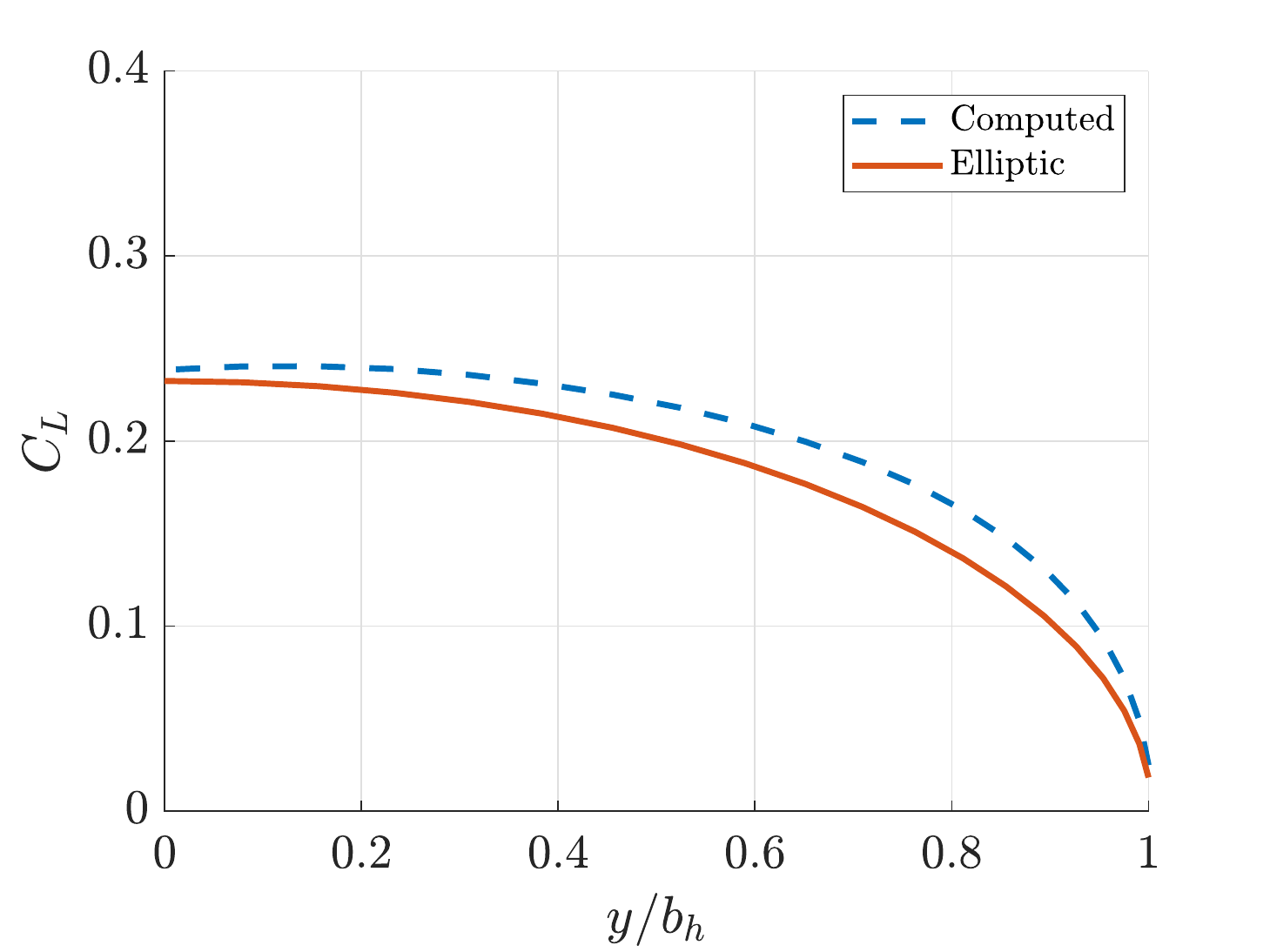}}
\caption{$C_L$ distribution of design 2}
\label{fig:CLdist_des2}
\end{figure}

\begin{figure}[hbt!]
\centering
\includegraphics[width=0.49\linewidth]{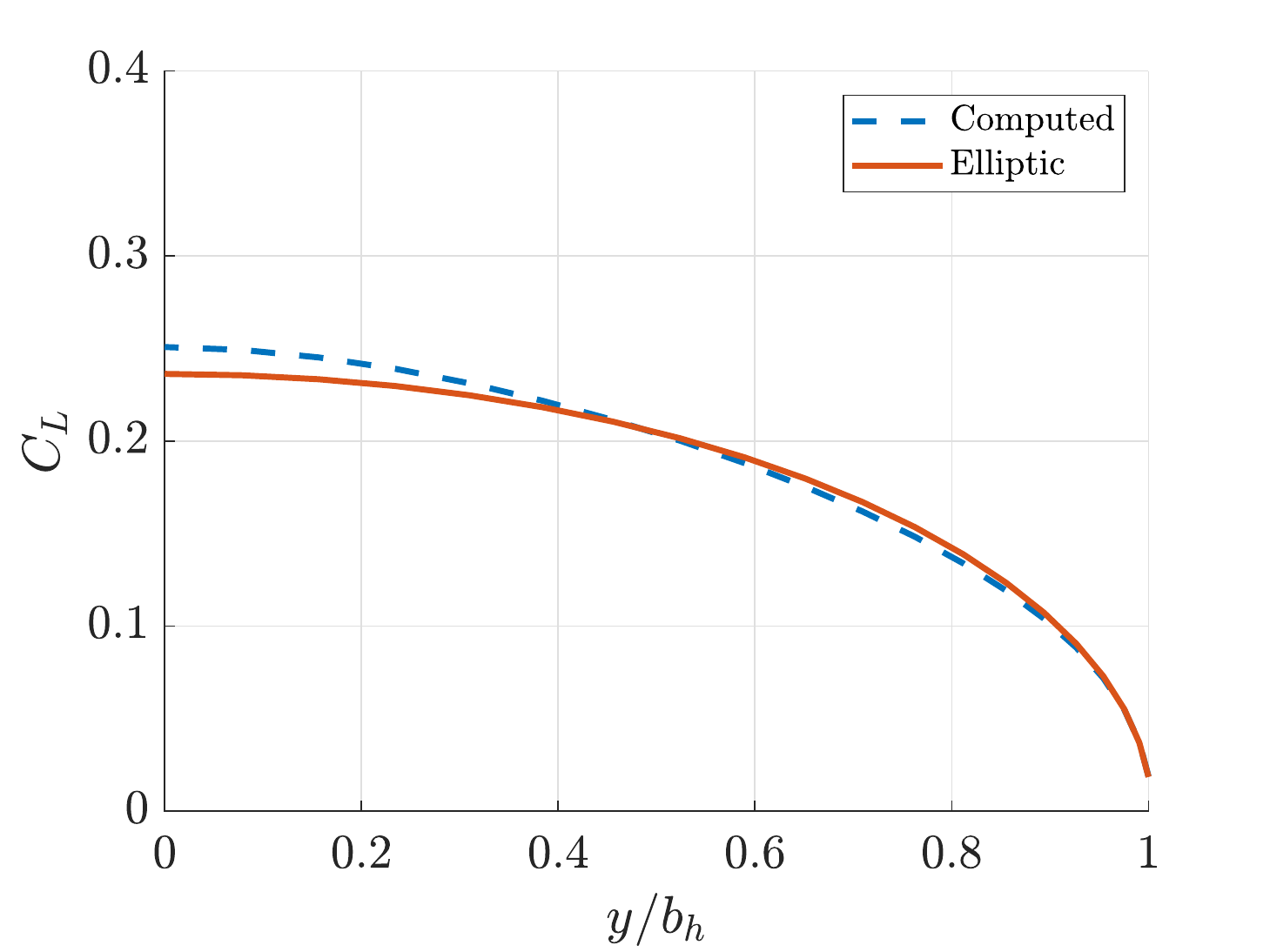}
\caption{Optimized $C_L$ distribution of design 1}
\label{fig:CLdist_des1}
\end{figure}

The lifting line theory can be theoretically employed to show that for minimum induced drag for a specific configuration, the required lift distribution for the wing should be elliptic in nature. Thus, such lift distribution is an indicator of the theoretical optimality of the configuration. We calculate the wing lift distribution for the initial and final configurations for design 2 to compare the effects of the DL-ROM-based optimization. The computed lift distribution compared to an elliptic lift distribution is shown in Figs.~\ref{fig:CLdist_des2} (a) and (b), for the initial and optimized designs, respectively. We can see that the initial design 2 not only violates certain constraints, but its lift distribution is far away from the optimal elliptic distribution, especially at the wing root. For the final design obtained via the DL-ROM-based optimization, the lift distribution does not exactly match the elliptic lift distribution but is much closer to it as a result of the optimization. To further understand the source of this difference in the optimized lift distribution, we perform a similar comparison for the optimized design 1. Since the optimized design 1 has a lower $C_{D_i}$ than optimized design 2, we can expect a better match with the theoretical prediction for a global optimum. Fig~\ref{fig:CLdist_des1} shows that the lift distribution of the optimized design 1 matches the elliptic lift distribution almost exactly near the wingtip. Some small differences are observed near the wing root. Overall, it follows the theoretically predicted optimum much better than design 2. This shows that DL-ROM-based optimization can provide physically consistent results over a range of design values. 

The small differences between the theoretical and actual optimized solutions for design 1 could be a consequence of several factors. It could arise from the assumptions considered in the development of the DL-ROM ranging from implicit shape approximation via level-sets to the FOM approximation via CNNs. However, it is also possible that the exploration of the design space with gradient-based optimization could not lead to the true global minimum. Thus, future shape optimization research would consider the application of metaheuristic global optimization algorithms \cite{samanipour2020adaptive} or a multi-fidelity analysis with trust region-based algorithm \cite{marques2021non}, both of which offer more robustness than gradient-based methods. 

\subsection{Computational efficiency}
The total cost of training the CNN model was about 11 hours on a single NVIDIA v100 Tensor Core GPU. The prediction phase for the CNN model includes the generation of the distance function from the airfoil shape parameters and subsequent CNN-based computation. This requires approximately 30 CPU seconds. The SQP-based gradient optimization with the CNN-based DL-ROM required roughly 500 CPU seconds on average. Since the CNN-based optimization is a completely data-driven process, its total computation time is independent of the fidelity of the computational solver or the experimental process employed. The reported computational time associated with two-dimensional RANS-based optimization \cite{bouhlel2020scalable} and direct numerical simulations of two-dimensional Navier-Stokes equations \cite{miyanawala2017efficient} indicate that the DL-ROM can provide significant computational efficiency during online applications even for two-dimensional cases. For the three-dimensional cases shown here, the reduction in computational time can be expected to increase many folds, especially as the number of design points increases. 

\section{Conclusions}
In this paper, we presented a novel deep learning-based model order reduction methodology for shape optimization applications via a deep convolutional neural network (CNN) and a level-set shape representation. The CNN provides a nonlinear mapping between the shapes and their associated features while conserving the equivariance of these features to the shape translations. Such a mapping can be achieved via the CNN with sufficient accuracy but at an improved computational efficiency compared to high-fidelity full-order models. To provide an implicit representation of complex three-dimensional shapes on a uniform structured grid, level-set methods were employed. The level-set representation not only enables the generic application of CNN models but also provides coarse-graining when the explicit representation of the shapes requires many parameters and is not computationally efficient for optimization applications.

To show the efficacy of our proposed DL-ROM-based optimization methodology, we considered the problem of minimizing the induced drag of a wing in three-dimensional potential flow. We illustrated satisfactory training, validation and test performance of our deep CNN model by optimal hyperparameter tuning. We also showed a reasonably accurate comparison between the aerodynamic coefficient of our optimized solutions to those obtained from the potential flow solver used for generating the three-dimensional training data sets. Furthermore, we explored  the physically consistent nature of our DL-ROM-based optimization by demonstrating a close match between the optimal design lift distribution obtained via the DL-ROM to the elliptic lift distribution predicted for a theoretical optimum. We also provided interpretability to the learning mechanism of the deep CNN model by exploring the kernel feature maps of the various convolutional layers of our deep CNN model. 

The proposed deep learning-based ROM can be applied non-intrusively during optimization, thus enabling seamless integration to general purpose gradient-based or metaheuristic optimization algorithms. 
The CNN-based shape optimization algorithm demonstrates computational efficiency during online optimization applications and promises a tractable solution for performing shape optimization for any general three-dimensional shape optimization problem. In the future, we plan to employ the DL-ROM-based optimization methodology to more physically complex and larger-scale problems, where training data will be obtained from high-fidelity computational solvers or experiments. We also plan to integrate our deep learning-based ROM as part of a multi-fidelity analysis model in the future with various gradient-free optimization techniques.

\section*{Acknowledgments}
This research was supported by the Natural Sciences and Engineering Research Council of Canada (NSERC) [grant number IRCPJ 550069-19]. The second author would like to thank the financial support from the NSERC Undergraduate Student Research Award. We would also like to acknowledge that the GPU facilities at the Compute Canada clusters were used for the training of our deep convolutional neural network models.

%
\section*{Conflict of interest}
The authors declare that they have no conflict of interest.

\section*{Replication of results}
Trained CNN models for the lift and induced drag predictions, CNN Matlab training scripts and input-output test data sets are available at \url{https://github.com/Wrik/CNN-optimization}.


\bibliographystyle{spmpsci}      
\bibliography{DL_optimization}

\begin{thebibliography}{10}
\providecommand{\url}[1]{{#1}}
\providecommand{\urlprefix}{URL }
\expandafter\ifx\csname urlstyle\endcsname\relax
  \providecommand{\doi}[1]{DOI~\discretionary{}{}{}#1}\else
  \providecommand{\doi}{DOI~\discretionary{}{}{}\begingroup
  \urlstyle{rm}\Url}\fi

\bibitem{allaire2004structural}
Allaire, G., Jouve, F., Toader, A.M.: Structural optimization using sensitivity
  analysis and a level-set method.
\newblock Journal of computational physics \textbf{194}(1), 363--393 (2004)

\bibitem{anderson2010fundamentals}
Anderson~Jr, J.D.: Fundamentals of aerodynamics.
\newblock Tata McGraw-Hill Education (2010)

\bibitem{blatman2009adaptive}
Blatman, G.: Adaptive sparse polynomial chaos expansions for uncertainty
  propagation and sensitivity analysis.
\newblock Ph.D. thesis, Clermont-Ferrand 2 (2009)

\bibitem{bonfiglio2018improving}
Bonfiglio, L., Perdikaris, P., Vernengo, G., de~Medeiros, J.S., Karniadakis,
  G.: Improving swath seakeeping performance using multi-fidelity gaussian
  process and bayesian optimization.
\newblock Journal of Ship Research \textbf{62}(04), 223--240 (2018)

\bibitem{bouhlel2016improving}
Bouhlel, M.A., Bartoli, N., Otsmane, A., Morlier, J.: Improving kriging
  surrogates of high-dimensional design models by partial least squares
  dimension reduction.
\newblock Structural and Multidisciplinary Optimization \textbf{53}(5),
  935--952 (2016)

\bibitem{bouhlel2020scalable}
Bouhlel, M.A., He, S., Martins, J.R.: Scalable gradient-enhanced artificial
  neural networks for airfoil shape design in the subsonic and transonic
  regimes.
\newblock Structural and Multidisciplinary Optimization pp. 1--14 (2020)

\bibitem{bronstein2017geometric}
Bronstein, M.M., Bruna, J., LeCun, Y., Szlam, A., Vandergheynst, P.: Geometric
  deep learning: going beyond euclidean data.
\newblock IEEE Signal Processing Magazine \textbf{34}(4), 18--42 (2017)

\bibitem{chernukhin2013multimodality}
Chernukhin, O., Zingg, D.W.: Multimodality and global optimization in
  aerodynamic design.
\newblock AIAA journal \textbf{51}(6), 1342--1354 (2013)

\bibitem{cybenko1989approximation}
Cybenko, G.: Approximation by superpositions of a sigmoidal function.
\newblock Mathematics of control, signals and systems \textbf{2}(4), 303--314
  (1989)

\bibitem{de2019structural}
De, S., Jrad, M., Kapania, R.K.: Structural optimization of internal structure
  of aircraft wings with curvilinear spars and ribs.
\newblock Journal of Aircraft \textbf{56}(2), 707--718 (2019)

\bibitem{de2021lightweight}
De, S., Singh, K., Seo, J., Kapania, R.K., Ostergaard, E., Angelini, N.,
  Aguero, R.: Lightweight chassis design of hybrid trucks considering multiple
  road conditions and constraints.
\newblock World Electric Vehicle Journal \textbf{12}(1), 3 (2021)

\bibitem{fan2019reliability}
Fan, X., Wang, P., Hao, F.: Reliability-based design optimization of crane
  bridges using kriging-based surrogate models.
\newblock Structural and Multidisciplinary Optimization \textbf{59}(3),
  993--1005 (2019)

\bibitem{feng2018multidisciplinary}
Feng, Y., Chen, Z., Dai, Y., Wang, F., Cai, J., Shen, Z.: Multidisciplinary
  optimization of an offshore aquaculture vessel hull form based on the support
  vector regression surrogate model.
\newblock Ocean Engineering \textbf{166}, 145--158 (2018)

\bibitem{garg2015high}
Garg, N., Kenway, G.K., Lyu, Z., Martins, J.R., Young, Y.L.: High-fidelity
  hydrodynamic shape optimization of a 3-d hydrofoil.
\newblock Journal of Ship Research \textbf{59}(4), 209--226 (2015)

\bibitem{Giles_2000}
Giles, M.B., Pierce, N.A.: An introduction to the adjoint approach to design.
\newblock Flow, Turbulence and Combustion \textbf{65}, 393–415 (2000)

\bibitem{goodfellow2016deep}
Goodfellow, I., Bengio, Y., Courville, A.: Deep learning.
\newblock MIT press (2016)

\bibitem{haftka1986structural}
Haftka, R.T., Grandhi, R.V.: Structural shape optimization—a survey.
\newblock Computer methods in applied mechanics and engineering \textbf{57}(1),
  91--106 (1986)

\bibitem{han2015surrogate}
Han, Z.H., Abu-Zurayk, M., G{\"o}rtz, S., Ilic, C.: Surrogate-based aerodynamic
  shape optimization of a wing-body transport aircraft configuration.
\newblock In: Symposium on AeroStructures, pp. 257--282. Springer (2015)

\bibitem{he2019robust}
He, X., Li, J., Mader, C.A., Yildirim, A., Martins, J.R.: Robust aerodynamic
  shape optimization—from a circle to an airfoil.
\newblock Aerospace Science and Technology \textbf{87}, 48--61 (2019)

\bibitem{hornik1991approximation}
Hornik, K.: Approximation capabilities of multilayer feedforward networks.
\newblock Neural networks \textbf{4}(2), 251--257 (1991)

\bibitem{ioffe2015batch}
Ioffe, S., Szegedy, C.: Batch normalization: Accelerating deep network training
  by reducing internal covariate shift.
\newblock In: International conference on machine learning, pp. 448--456. PMLR
  (2015)

\bibitem{jameson_1998}
Jameson, A., Pierce, N., Martinelli., L.: Optimum aerodynamic design using the
  navier-stokes equations.
\newblock J. Theor. Comp. Fluid Mech \textbf{10}, 213--237 (1998)

\bibitem{kapania2018multidisciplinary}
Kapania, R.K., Schetz, J.A., Mallik, W., Segee, M.C., Gupta, R.:
  Multidisciplinary design optimization and cruise mach number study of
  truss-braced wing aircraft.
\newblock Tech. rep., NASA/CR-2018-219836, NF1676L-27943 (2018)

\bibitem{keane2008engineering}
Keane, A., Forrester, A., Sobester, A.: Engineering design via surrogate
  modelling: a practical guide.
\newblock American Institute of Aeronautics and Astronautics, Inc. (2008)

\bibitem{kingma2014adam}
Kingma, D.P., Ba, J.: Adam: A method for stochastic optimization.
\newblock arXiv preprint arXiv:1412.6980  (2014)

\bibitem{koziel2013surrogate}
Koziel, S., Leifsson, L.: Surrogate-based aerodynamic shape optimization by
  variable-resolution models.
\newblock AIAA journal \textbf{51}(1), 94--106 (2013)

\bibitem{krige1951statistical}
Krige, D.G.: A statistical approach to some basic mine valuation problems on
  the witwatersrand.
\newblock Journal of the Southern African Institute of Mining and Metallurgy
  \textbf{52}(6), 119--139 (1951)

\bibitem{lecun2015deep}
LeCun, Y., Bengio, Y., Hinton, G.: Deep learning.
\newblock Nature \textbf{521}(7553), 436--444 (2015)

\bibitem{lyu2015aerodynamic}
Lyu, Z., Kenway, G.K., Martins, J.R.: Aerodynamic shape optimization
  investigations of the common research model wing benchmark.
\newblock AIAA journal \textbf{53}(4), 968--985 (2015)

\bibitem{maas2013rectifier}
Maas, A.L., Hannun, A.Y., Ng, A.Y.: Rectifier nonlinearities improve neural
  network acoustic models.
\newblock In: Proc. icml, vol.~30, p.~3. Citeseer (2013)

\bibitem{mallik2015effect}
Mallik, W., Kapania, R.K., Schetz, J.A.: Effect of flutter on the
  multidisciplinary design optimization of truss-braced-wing aircraft.
\newblock Journal of Aircraft \textbf{52}(6), 1858--1872 (2015)

\bibitem{mallik2020kriging}
Mallik, W., Raveh, D.E.: Kriging-based aeroelastic gust response analysis at
  high angles of attack.
\newblock AIAA Journal \textbf{58}(9), 3777--3787 (2020)

\bibitem{marques2021non}
Marques, S.P., Kob, L., Robinson, T.T., Yao, W., Sun, L.: Non-intrusive
  aerodynamic shape optimisation with a discrete empirical interpolation
  method.
\newblock In: AIAA Scitech 2021 Forum, p. 0172 (2021)

\bibitem{martins2013review}
Martins, J.R., Hwang, J.T.: Review and unification of methods for computing
  derivatives of multidisciplinary computational models.
\newblock AIAA journal \textbf{51}(11), 2582--2599 (2013)

\bibitem{matthews2018gaussian}
Matthews, A.G.d.G., Rowland, M., Hron, J., Turner, R.E., Ghahramani, Z.:
  Gaussian process behaviour in wide deep neural networks.
\newblock arXiv preprint arXiv:1804.11271  (2018)

\bibitem{miyanawala2018novel}
Miyanawala, T., Jaiman, R.K.: A novel deep learning method for the predictions
  of current forces on bluff bodies.
\newblock In: International Conference on Offshore Mechanics and Arctic
  Engineering, vol. 51210, p. V002T08A003. American Society of Mechanical
  Engineers (2018)

\bibitem{miyanawala2017efficient}
Miyanawala, T.P., Jaiman, R.K.: An efficient deep learning technique for the
  navier-stokes equations: Application to unsteady wake flow dynamics.
\newblock arXiv preprint arXiv:1710.09099  (2017)

\bibitem{pironneau_2004}
Mohammadi, B., Pironneau, O.: Shape optimization in fluid mechanics.
\newblock Annual Review of Fluid Mechanics \textbf{36}(1), 255--279 (2004)

\bibitem{mojgani2021low}
Mojgani, R., Balajewicz, M.: Low-rank registration based manifolds for
  convection-dominated pdes.
\newblock In: Proceedings of the AAAI Conference on Artificial Intelligence,
  vol.~35, pp. 399--407 (2021)

\bibitem{nagawkar2022aerodynamic}
Nagawkar, J.R., Leifsson, L.T., He, P.: Aerodynamic shape optimization using
  gradient-enhanced multifidelity neural networks.
\newblock In: AIAA SCITECH 2022 Forum, p. 2350 (2022)

\bibitem{perdikaris2017nonlinear}
Perdikaris, P., Raissi, M., Damianou, A., Lawrence, N.D., Karniadakis, G.E.:
  Nonlinear information fusion algorithms for data-efficient multi-fidelity
  modelling.
\newblock Proceedings of the Royal Society A: Mathematical, Physical and
  Engineering Sciences \textbf{473}(2198), 20160751 (2017)

\bibitem{raissi2016deep}
Raissi, M., Karniadakis, G.: Deep multi-fidelity gaussian processes.
\newblock arXiv preprint arXiv:1604.07484  (2016)

\bibitem{reuther1996aerodynamic}
Reuther, J., Jameson, A., Farmer, J., Martinelli, L., Saunders, D.: Aerodynamic
  shape optimization of complex aircraft configurations via an adjoint
  formulation.
\newblock In: 34th Aerospace Sciences Meeting and Exhibit, p.~94 (1996)

\bibitem{saha2022practical}
Saha, D., De, S.: Practical self-driving cars: Survey of the state-of-the-art
  (2022)

\bibitem{samanipour2020adaptive}
Samanipour, F., Jelovica, J.: Adaptive repair method for constraint handling in
  multi-objective genetic algorithm based on relationship between constraints
  and variables.
\newblock Applied Soft Computing \textbf{90}, 106143 (2020)

\bibitem{samareh2004aerodynamic}
Samareh, J.: Aerodynamic shape optimization based on free-form deformation.
\newblock In: 10th AIAA/ISSMO multidisciplinary analysis and optimization
  conference, p. 4630 (2004)

\bibitem{schobi2015polynomial}
Schobi, R., Sudret, B., Wiart, J.: Polynomial-chaos-based kriging.
\newblock International Journal for Uncertainty Quantification \textbf{5}(2)
  (2015)

\bibitem{sethian1996theory}
Sethian, J.A.: Theory, algorithms, and applications of level set methods for
  propagating interfaces.
\newblock Acta numerica \textbf{5}, 309--395 (1996)

\bibitem{singh2021alga}
Singh, K., Kapania, R.K.: Alga: Active learning-based genetic algorithm for
  accelerating structural optimization.
\newblock AIAA Journal \textbf{59}(1), 330--344 (2021)

\bibitem{sobester2008quest}
S{\'o}bester, A., Barrett, T.: Quest for a truly parsimonious airfoil
  parameterization scheme.
\newblock In: The 26th Congress of ICAS and 8th AIAA ATIO, p. 8879 (2008)

\bibitem{sobester2014aircraft}
S{\'o}bester, A., Forrester, A.I.: Aircraft aerodynamic design: geometry and
  optimization.
\newblock John Wiley \& Sons (2014)

\bibitem{streuber2020evaluating}
Streuber, G.M., Zingg, D.W.: Evaluating the risk of local optima in aerodynamic
  shape optimization.
\newblock AIAA Journal pp. 1--13 (2020)

\bibitem{matlabdl}
The~MathWorks, I.: Matlab deep learning toolbox.
\newblock \urlprefix\url{https://www.mathworks.com/help/deeplearning/}

\bibitem{matlabopt}
The~MathWorks, I.: Matlab optimization toolbox.
\newblock \urlprefix\url{https://www.mathworks.com/help/optim/}

\bibitem{timme2011transonic}
Timme, S., Marques, S., Badcock, K.: Transonic aeroelastic stability analysis
  using a kriging-based schur complement formulation.
\newblock AIAA journal \textbf{49}(6), 1202--1213 (2011)

\bibitem{wang2007extended}
Wang, S., Lim, K.M., Khoo, B.C., Wang, M.Y.: An extended level set method for
  shape and topology optimization.
\newblock Journal of Computational Physics \textbf{221}(1), 395--421 (2007)

\bibitem{wang2018database}
Wang, S., Sun, G., Chen, W., Zhong, Y.: Database self-expansion based on
  artificial neural network: An approach in aircraft design.
\newblock Aerospace Science and Technology \textbf{72}, 77--83 (2018)

\bibitem{wang2006radial}
Wang, S., Wang, M.Y.: Radial basis functions and level set method for
  structural topology optimization.
\newblock International journal for numerical methods in engineering
  \textbf{65}(12), 2060--2090 (2006)

\bibitem{yao2020reduced}
Yao, W., Marques, S., Robinson, T., Armstrong, C., Sun, L.: A reduced-order
  model for gradient-based aerodynamic shape optimisation.
\newblock Aerospace Science and Technology \textbf{106}, 106120 (2020)

\bibitem{yondo2018review}
Yondo, R., Andr{\'e}s, E., Valero, E.: A review on design of experiments and
  surrogate models in aircraft real-time and many-query aerodynamic analyses.
\newblock Progress in aerospace sciences \textbf{96}, 23--61 (2018)

\end{thebibliography}

\end{document}